\documentclass[12pt]{iopart}

\usepackage[margin=2.5cm]{geometry}
\usepackage{setspace}
\usepackage{lmodern}
\expandafter\let\csname equation*\endcsname\relax
\expandafter\let\csname endequation*\endcsname\relax
\usepackage{amsmath,amssymb}
 \usepackage{amsfonts,amsthm}
 \usepackage{amscd,amsxtra,latexsym}
\usepackage{graphicx}
\usepackage{subfigure}
\usepackage{epsf,epsfig}
\usepackage{bm,bbm}
\usepackage{algpseudocode}
\usepackage{algorithmicx, algorithm}
\usepackage{physics}
\usepackage{color}
\usepackage{cite}
\usepackage{pxfonts}
\usepackage[footnotesize,bf]{caption}
\usepackage{hhline}
\usepackage[T1]{fontenc}
\usepackage{placeins}
\usepackage{multirow,dcolumn}
\usepackage{mathrsfs}
\usepackage{verbatim}
\usepackage{comment}
\usepackage{epstopdf}
\usepackage{lineno}
\usepackage[colorlinks=true,
linkcolor=blue,
anchorcolor=blue,
citecolor=blue,
urlcolor=blue,]{hyperref}
\newtheorem{remark}{Remark}[section]
\usepackage{iopams}
\begin{document}

\title[PI-INN]{Efficient Bayesian inference using physics-informed invertible neural networks for inverse problems}

\author{Xiaofei Guan\textsuperscript{1,2}, Xintong Wang\textsuperscript{1}, Hao Wu\textsuperscript{3,*}, Zihao Yang\textsuperscript{4}, Peng Yu\textsuperscript{5}}
\address{$^1$ School of Mathematical Sciences, Tongji University, Shanghai 200092, China}
\address{$^2$ Key Laboratory of Advanced Civil Engineering Materials of Ministry of Education, Tongji University, Shanghai 200092, China}
\address{$^3$ School of Mathematical Sciences, Institute of Natural Sciences, and MOE-LSC, Shanghai Jiaotong University, Shanghai 200240, China}
\address{$^4$ School of Mathematics and Statistics, Northwestern Polytechnical University, Xi'an 710129, China}
\address{$^5$ State Key Laboratory of Marine Geology, Tongji University, Shanghai 200092, China}
\ead{hwu81@sjtu.edu.cn}
%
%\address{IOP Publishing, Temple Circus, Temple Way, Bristol BS1 6HG, UK}
%\ead{submissions@iop.org}
\vspace{10pt}
\begin{indented}
\item[]September 2023
\end{indented}

\begin{abstract}
In this paper, we introduce an innovative approach for addressing Bayesian inverse problems through the utilization of physics-informed invertible neural networks (PI-INN). The PI-INN framework encompasses two sub-networks: an invertible neural network (INN) and a neural basis network (NB-Net). The primary role of the NB-Net lies in modeling the spatial basis functions characterizing the solution to the forward problem dictated by the underlying partial differential equation. Simultaneously, the INN is designed to partition the parameter vector linked to the input physical field into two distinct components: the expansion coefficients representing the forward problem solution and the Gaussian latent noise. If the forward mapping is precisely estimated, and the statistical independence between expansion coefficients and latent noise is well-maintained, the PI-INN offers a precise and efficient generative model for Bayesian inverse problems, yielding tractable posterior density estimates. As a particular physics-informed deep learning model, the primary training challenge for PI-INN centers on enforcing the independence constraint, which we tackle by introducing a novel independence loss based on estimated density. We support the efficacy and precision of the proposed PI-INN through a series of numerical experiments, including inverse kinematics, 1-dimensional and 2-dimensional diffusion equations, and seismic traveltime tomography. Specifically, our experimental results showcase the superior performance of the proposed independence loss in comparison to the commonly used but computationally demanding kernel-based maximum mean discrepancy loss.
\end{abstract}
%
% Uncomment for keywords
%\vspace{2pc}

%
% Uncomment for Submitted to journal title message
%\submitto{\JPA}
%
% Uncomment if a separate title page is required
%\maketitle
% 
% For two-column output uncomment the next line and choose [10pt] rather than [12pt] in the \documentclass declaration
%\ioptwocol

\section{Introduction}\label{section:intro}
Inverse problems arise in various scientific and engineering fields, such as seismic tomography \cite{Seismic-tomo}, material inspection \cite{Inspection-IP}, and medical imaging \cite{tomo-MI}, etc. These problems aim to infer the input parameters of a system, such as physical fields, source locations, initial and boundary conditions, from the corresponding measurement data. Inverse problems are typically ill-posed and difficult to solve because the measurement data originates from indirect and sparse observations. It is of great practical value and theoretical significance to develop an efficient computation model for accurate and stable solutions, which has drawn much attention from scientists and engineers.

Classical approaches to inverse problems encompass regularization methods and Bayesian inference. Regularization techniques aim to determine point estimates by simultaneously minimizing the misfits between observed and predicted outputs while penalizing undesired parameter features through a regularization term \cite{engl1996regularization,SampledTikhonov}. Bayesian inference incorporates uncertainties in prior information and observations through Bayesian rules. It constructs the posterior distribution of parameters, thus effectively quantifying the uncertainty of unknown inputs. However, in many applications, it is computationally expensive to sample from the intractable Bayesian posterior distributions by using statistical simulation techniques, e.g., Markov Chain Monte Carlo (MCMC)\cite{MultiscaleMRM, DataDriven,InvStekloff,IP-NU}, Hamiltonian Monte Carlo \cite{HMC-tomo,DD-HMC}, and Sequential Monte Carlo \cite{Multicale-SMC,SMC}. To address this issue, many researchers utilized various variational inference algorithms to improve computational efficiency, including mean-field variational inference \cite{ApproximateMI,BVI-Flow}, automatic differential variational inference\cite{VI-tomo,ADVI-Geo}, and Stein variational gradient descent\cite{VI-tomo,SVGD-inv}. But the variational methods may face challenges in high-dimensional settings due to the difficulty of accurately approximating posterior distributions with the tractable surrogate distribution. Although all aforementioned methods have shown to be relatively efficient, they still require numerous evaluations of the forward model and the complicated parametric derivatives, which introduce higher computational cost in high-dimensional situations. Consequently, deep learning-based methods offer an efficient alternative, where they are capable of providing real-time inversions for certain classes of inverse problems with new measurement data of similar type. 

Deep learning methods have emerged as a promising approach for solving inverse problems, such as medical imaging \cite{CNN-IP-image,Deep-MIA} and electromagnetic inversion \cite{Deep-EI,Deep-EI-1d}. Nevertheless, these methods heavily rely on labeled data from solutions of the forward problem, rendering them inappropriate for inverse problems where such information is not present.
With the continuous advancement of physics-informed neural network (PINN) technology in recent years, which leverages physical equations instead of labeled data to solve partial differential equations (PDEs) \cite{PINN}, extensive related research has demonstrated that PINN can efficiently recover unknown parameters as well \cite{Quantify-PINN,PINN-solid,DeepXDE,B-PINN,Adver-PINN,PINNtomo,PINN-wavefield,A-PINN,Surrogate-PINN}. Furthermore, models based on invertible neural networks (INN) \cite{NICE, RNVP, Glow} have exhibited potential in Bayesian inference. This is attributed to their efficiency and accuracy advantages in both sampling and density estimation of complex target distributions through bijective transformations, as evidenced in the works of \cite{INN, INN-geo, cINN, BayesFlow, NFFs, Inv-DeepONets, L-HYDRA, VI-NFs}.
Nonetheless, certain limitations persist within the current research.
For instance, INN models in \cite{INN, INN-geo, cINN, BayesFlow} demand substantial quantities of labeled data for pretraining, which is often unavailable in many practical inverse problems. The invertible DeepONet introduced in \cite{Inv-DeepONets} alleviates the need for labeled data, but it entails an additional variational inference procedure to approximate the true posterior distribution.

In this paper, we propose a novel approach to address the issues of Bayesian inverse problems, namely physics-informed invertible neural networks (PI-INN). In PI-INN, a neural Basis Network (NB-Net) identifies spatial basis functions for the solution to the forward problem governed by the given PDE, while an INN bijectively transforms the input physical field's parameter vector into a combination of expansion coefficients for the forward problem solution and Gaussian latent variables. To achieve an effective approximation of the true posterior, PI-INN undergoes training to accurately predict the solution of the equation while ensuring statistical independence between the expansion coefficients and Gaussian latent variables. Notably, we introduce a new independence loss term that can be naturally incorporated with residual based physics-informed loss terms, and theoretically prove that the PI-INN can efficiently generate samples from the true posterior distribution of unknown parameters at the limit case where all loss terms approaches zero.
The proposed PI-INN approach admits the following main advantages: 

\begin{itemize}
%	\item Our theoretical proof shows that PI-INN can efficiently generate samples from the true posterior distribution of the inverse problem when the loss function approaches zero. 
	\item PI-INN can provide a tractable estimation of the posterior distribution, which enables efficient sampling and high-accuracy density evaluation without requiring additional variational inference procedures \cite{Inv-DeepONets}.
	\item The innovative independence loss term leverages the INN's density estimation capabilities, obviating the requirement for labeled data. In contrast, the kernel-based MMD loss term employed in \cite{INN, INN-geo}  is applicable exclusively to labeled data and easily suffers from the curse of dimensionality.
	\item PI-INN establishes an integrated learning framework for Bayesian inverse problems, with the capability to merge physics-informed learning and data-driven learning.
\end{itemize}

The remainder of this paper is organized as follows. Section \ref{section:definition} provides the problem statement of Bayesian inverse problems. In Section \ref{section:methody}, the framework of our physics-informed invertible neural networks (PI-INN) is introduced, where the algorithm of training and inversion approaches with PI-INN for the physical system is described in detail. In section \ref{section:experiments}, we first present the advantage of the new loss term and then demonstrate the performance of PI-INN for solving Bayesian inverse problems through several numerical examples. Finally, some concluding remarks are given in Section \ref{section:conclusions}.

\section{Problem statement}\label{section:definition}
Consider a physical system described by the following PDE:
\begin{equation}\label{general-equ}
\left\{
\begin{array}{ll}
\mathcal{N}(u(\bm{x});K(\bm{x};\lambda)) = 0,\quad \bm{x}\in\mathcal{D},\\
\mathcal{B}(u(\bm{x})) = 0,\quad \bm{x}\in\partial\mathcal{D},
\end{array}
\right.
\end{equation}
where $\mathcal{N}$ denotes the general differential operator defined in the domain $\mathcal{D}\subset\mathbb{R}^{d}\ (d = 1,2, $\text{or} 3$)$, $\mathcal{B}$ is the boundary operator on the boundary $\partial\mathcal{D}$, $K(\bm{x};\lambda)$ represents the input physical field characterized by a vector $\lambda$ and $u(\bm{x})$ denotes the solutions of the PDE. 

\begin{figure}[htpb]
	\centering
	\includegraphics[scale=0.45]{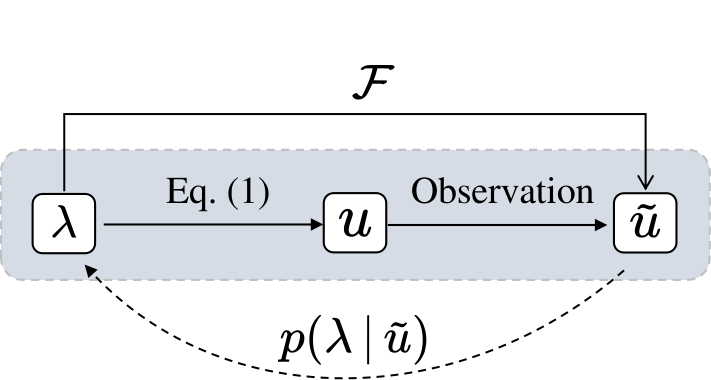}
	\caption{Graphical representation of the Bayesian inversion method.}
	\label{graphical-setting}
\end{figure}

The corresponding inverse problem aims to recover $\lambda$ from the observations $\tilde{u}\in\mathbb{R}^{D}$ of $u$ at specific sensors of the domain. The relationship of the above variables in observation system is illustrated in Fig. \ref{graphical-setting}. For the Bayesian inference method, the target is to obtain samples from the posterior distribution $p(\lambda\,|\,\tilde{u})$:
\begin{equation}\label{posterior}
p(\lambda\,|\,\tilde{u})\propto p(\tilde{u}\,|\,\lambda)\,p(\lambda),
\end{equation}
where $p(\lambda)$ represents the prior and $p(\tilde{u}\,|\,\lambda)$ is the likelihood. 
The main challenge associated with the classical Bayesian approaches lies in the computational burden of repeatly evaluating the forward model $\mathcal{F}$ for Monte Carlo sampling. To tackle this challenge, we introduce an INN-based model that enables accurate and efficient posterior sampling, bypassing the need for explicit application of Bayes' rule.

\section{Methodology}\label{section:methody}
\subsection{Invertible neural networks}
\label{INN-NF}
\begin{figure}[htp]
	\centering
	\subfigure[]{
		\includegraphics[scale=0.69]{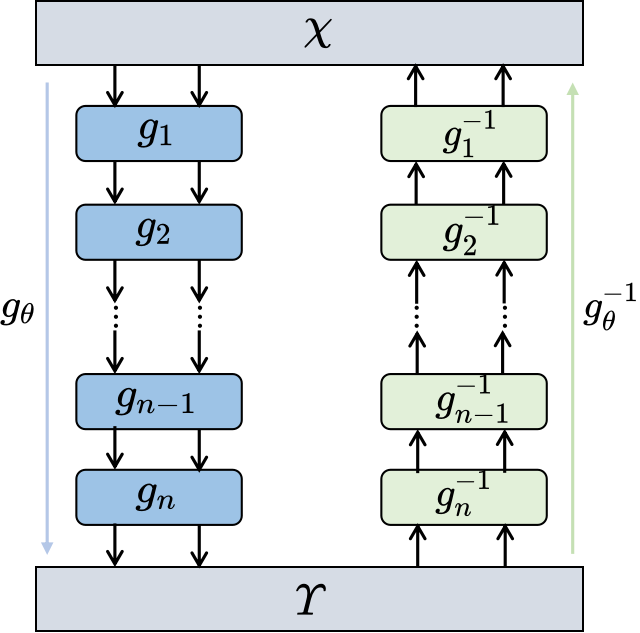}
	}
	\hspace{+25pt}
	\subfigure[]{
		\includegraphics[scale=0.53]{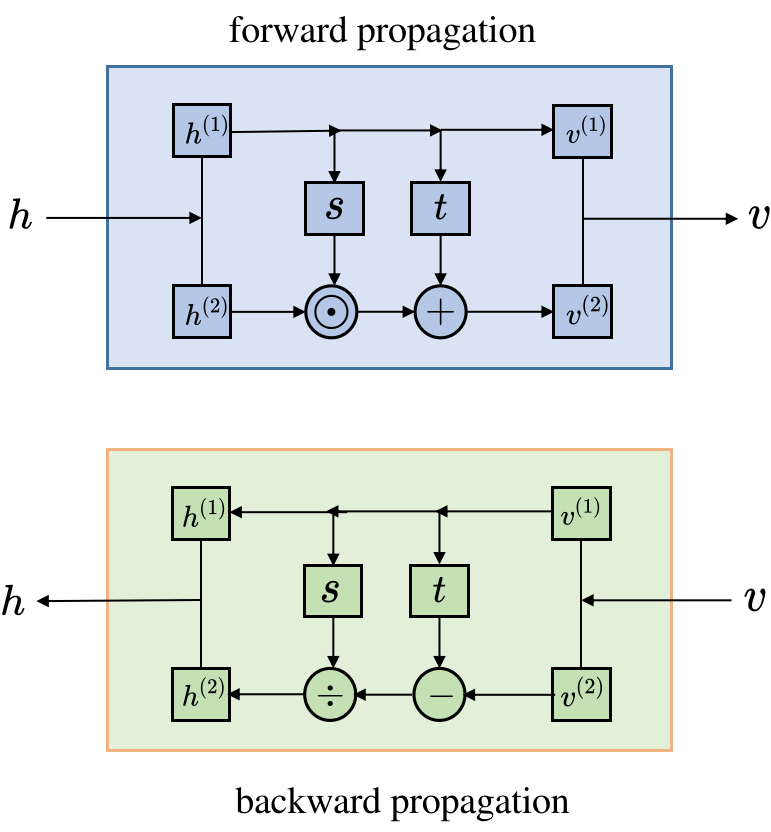}
	}
	\caption{(a): INN with $n$ coupling layers; (b): Illustration of affine coupling layers.}
	\label{INNs}
\end{figure}

Invertible neural networks (INNs) \cite{NICE, RNVP, Glow}, also known as flow-based models, are a particular type of neural network that can build bijective transformations between inputs and outputs. As shown in Fig. \ref{INNs}\,(a), INNs build the bijection $g_{\bm{\theta}}$ by stacking a series of reversible and differentiable coupling layers $g_1,\cdots,g_n$, where $\varUpsilon\in\mathbb{R}^{F}$ and $\chi\in\mathbb{R}^{F}$ are the random variables, $\varUpsilon = g_{\theta}(\chi)$. The probability distributions of $\chi$ and $\varUpsilon$ satisfy
\begin{equation}\label{pdf-series-formula}
\mathrm{log}\,p_{\chi}(\chi) = \mathrm{log}\,p_{\varUpsilon}(g_{\bm{\theta}}(\chi)) - \sum\limits_{j=1}^{n}\mathrm{log}\,\Big\vert\mathrm{det}\dfrac{\mathrm{d} g_{j}(h_{j-1})}{\mathrm{d}h_{j-1}}\Big\vert,
\end{equation}
where $h_j$ is the intermediate variable such that $h_j = g_j(h_{j-1}), j=1,\cdots,n$, $h_0 = \chi$, $h_n = \varUpsilon$.

In this work, we utilize the affine coupling layer as the basic building block, which was introduced in \cite{RNVP, Glow}. Fig.\,\ref{INNs}\,(b) illustrates the forward and backward operations of the affine coupling layer. In the forward propagation, the input $h\in\mathbb{R}^{F}$ is first split into two parts: $h^{(1)}\in\mathbb{R}^{f}$ and $h^{(2)}\in\mathbb{R}^{F-f},\ f < F$. These two parts then undergo the coupling operation as follows:
\begin{equation}\label{AF-forward}
\centering
v^{(1)} = h^{(1)},\quad v^{(2)} = h^{(2)}\odot\mathrm{exp}(\bm{s}(h^{(1)})) + \bm{t}(h^{(2)}),
\end{equation}
where $\bm{s}(\cdot)$ and $\bm{t}(\cdot)$ are scale and translation functions respectively, and they are both modeled by trainable neural networks. Here, $\odot$ represents element-wise multiplication, and the inverse transformation can be trivially given:
\begin{equation}\label{AF-backward}
h^{(1)} = v^{(1)},\quad h^{(2)} = (v^{(2)} - \bm{t}(v^{(1)}))\odot\mathrm{exp}(-\bm{s}(v^{(1)})).
\end{equation}
The Jacobian of the forward map is given as
\begin{equation}\label{AF-jacobi}
\frac{\partial(v^{(1)},v^{(2)})}{\partial(h^{(1)},h^{(2)})} = \begin{bmatrix}
\mathbb{I}_{d} & \bm{0}_{f\times (F-f)}\\
\dfrac{\partial v^{(2)}}{\partial h^{(1)}} & \mathrm{diag}(\mathrm{exp}(\bm{s}(h^{(1)})))
\end{bmatrix}.
\end{equation}
Due to the tractable inverse map and determinant of the Jacobian for each layer, the change-of-variable formula (\ref{pdf-series-formula}) can be easily calculated.

%INNs possess the capability to learn the forward mapping of a physical system while concurrently addressing inverse problems as shown in Fig.~\ref{graphical-setting}, thanks to their inherent invertibility. Nevertheless, in real-world scenarios, inverse problems frequently lack unique solutions due to information loss from measurements, preventing the establishment of a bijection between measurements and input parameters.
%Additionally, the scarcity of labeled data presents challenges when learning the forward mapping.
%To address this challenge, we introduce physics-informed invertible neural networks (PI-INN) in the following section.
In the context of this study, a pertinent question arises: Can we harness the capabilities of INNs to acquire the forward mapping of a physical system while simultaneously tackling inverse problems, as illustrated in Fig.~\ref{graphical-setting}. Regrettably, practical real-world scenarios often entail inverse problems devoid of unique solutions, primarily due to the loss of information in measurements, which hinders the establishment of a bijection between measurements and input parameters. To confront this challenge, we introduce physics-informed invertible neural networks (PI-INN) in the subsequent section.

\subsection{Physics-informed invertible neural networks}
\label{networks}
\begin{figure}[htp]
	\centering
	\includegraphics[scale=0.55]{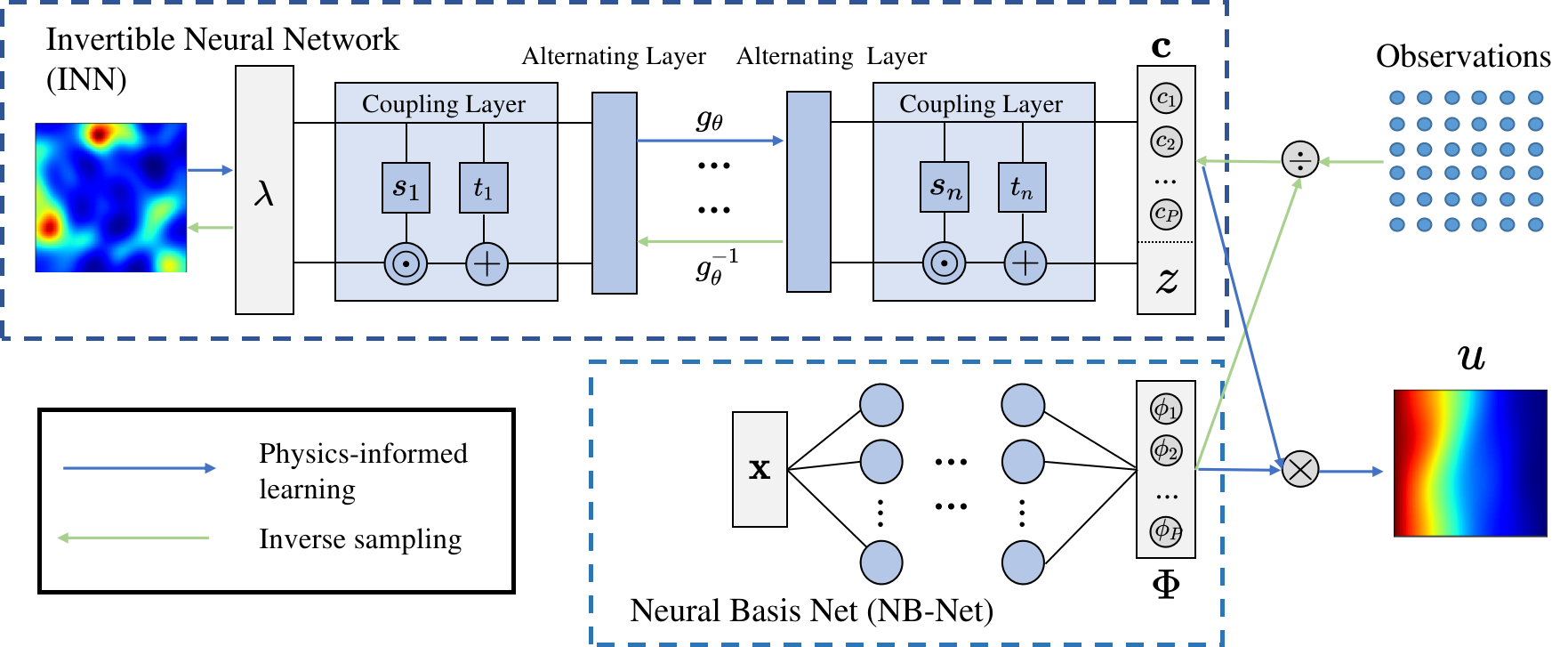}
	\caption{Schematic of PI-INN. The INN coupling with NB-Net can not only obtain the solutions of forward problems but also infer the unknown parameters from the observations. The blue and green solid lines represent the directions of learning and Bayesian inversion processes, respectively. The standard Gaussian distribution is chosen as the prior distribution for $z$ in this paper.}
	\label{PI-INN}
\end{figure}

The architecture of PI-INN is depicted in Fig.~\ref{PI-INN}, which consists of an INN and a neural basis net (NB-Net). Here, we consider the separate representation of the forward problem solution akin to the representation taken in DeepONets \cite{DeepONet, PI-DeepONet}:
\begin{equation}\label{solution-expansion}
u(\bm{x};\lambda) = \sum\limits_{i=1}^{P}c_i(\lambda)\phi_i(\bm{x}),
\end{equation}
where $\bm c=(c_1,\ldots,c_P)^\top$ represents the expansion coefficient, and $\phi_i(\bm{x})$ is the basis function defined on the spatial domain. 

The INN $g_\theta$, which consists of affine coupling layers, is designed to learn the bijective mapping between the input parameter $\lambda$ and the pair $[\bm{c}, z]$ of expansion coefficients and random latent variables. The random latent variable $z$ is employed to capture the inherent information that is not encompassed by $\bm{c}$, where $\mathrm{ndim}(z)=\mathrm{ndim}(\lambda) - \mathrm{ndim}(\bm{c})$. It is assumed that $z$ is statistically independent of $\bm{c}$ and follows a tractable probability distribution, e.g., standard Gaussian distribution. Consequently, the forward map from $\lambda$ to $\bm c$ dfined by $g_{\bm{\theta}}$ is deterministic as shown in Fig.~\ref{PI-INN}, whereas the inverse map contains uncertainty.

\begin{remark}
In this study, we make the assumption that $\mathrm{ndim}(\lambda) > \mathrm{ndim}(\bm{c})$, which aligns with the typical scenario in practical applications where the dimensionality of measurement variables is lower than that of the parameters. 
\end{remark}

The NB-Net is a fully-connected neural network with $P$ output neurons, utilized for modeling spatial basis functions $\phi=(\phi_1,\ldots,\phi_P)$ . For a given $\phi$, PI-INN is capable of predicting the numerical solutions at any positions of the domain $\mathcal{D}$ and solving inverse problems with variable locations of sensors. According to Eq.~(\ref{solution-expansion}), the vector representation of the solutions on $\{\bm{x}^{(i)}\}_{i=1}^{M}\subset\mathcal{D}$ is defined as follows:
\begin{equation}\label{solution-vector-representation}
\bm{u}=\Phi\bm{c},
\end{equation}
where
$$\Phi = \begin{bmatrix}
\phi_1(\bm{x}^{(1)}) & \cdots & \phi_P(\bm{x}^{(1)})\\
\vdots & \vdots & \vdots\\
\phi_1(\bm{x}^{(M)}) & \cdots & \phi_P(\bm{x}^{(M)})
\end{bmatrix},\quad \bm{c} = \begin{bmatrix}
c_1\\ \vdots \\ c_P
\end{bmatrix},\quad \bm{u} = \begin{bmatrix}
u(\bm{x}^{(1)})\\ \vdots \\ u(\bm{x}^{(M)})
\end{bmatrix}.$$
Here, $\{\bm{x}^{(i)}\}_{i=1}^{M}$ is an arbitrary spatital coordinates in $\mathcal{D}$, and $\Phi\in\mathbb{R}^{M\times P}$ can be calculated by the NB-Net. 

Meanwhile, given measurements $\tilde{\bm{u}} = [\tilde{u} (\tilde{\bm{x}}^{(1)}),\cdots,\tilde{u}(\tilde{\bm{x}}^{(O)})]^{\top}$ on $\{\tilde{\bm{x}}^{(i)}\}_{i=1}^{O}$, $\bm{c}$ can be deterimined by solving the following optimization problem \cite{L-HYDRA}:
\begin{equation}\label{lse-problem}
\mathop{\mathrm{min}}_{\bm{c}}\Vert\tilde{\bm{u}} - \tilde{\Phi}\bm{c}\Vert^2 - \rho\mathrm{log}\,p(\bm{c}),
\end{equation}
where $[\tilde{\Phi}]_{ij}=\phi_i(\tilde x^{(j)})$, $\rho$ is the regularization coefficient, $p(\bm{c})$ denotes the marginal distribution of $\bm{c}$, and can also be calculated by PI-INN (see below).

In summary, the forward and inverse implementation of PI-INN can be succinctly described as:
\[
\begin{array}{lcccccc}
\text{\textbf{Forward:}} & \quad & \lambda & \stackrel{g_{\theta}}{\xrightarrow[\quad]{\hspace*{1cm}}} & \left(\begin{array}{c}
\bm{c}\\
z
\end{array}\right) & \begin{array}{c}
\stackrel{\text{Eq. }(8)}{\xrightarrow[\quad]{\hspace*{1cm}}}\\
\\
\end{array} & \begin{array}{c}
\bm u\\
\\
\end{array}\\
\text{\textbf{Inverse:}} & \quad & \begin{array}{c}
\tilde{\bm u}\\
\mathcal{N}(0,I)
\end{array} & \begin{array}{c}
\stackrel{\text{Eq. }(9)}{\xrightarrow[\quad]{\hspace*{1cm}}}\\
\stackrel{\quad}{\xrightarrow[\quad]{\hspace*{1cm}}}
\end{array} & \left(\begin{array}{c}
\bm{c}\\
z
\end{array}\right) & \stackrel{g_{\theta}^{-1}}{\xrightarrow[\quad]{\hspace*{1cm}}} & \lambda
\end{array}
\]
The detailed training and inversion processes will be introduced by using PI-INN in the following section.

\subsection{Training and inversion procedure of PI-INN}
\label{training-test-approach}

During the learning process of PI-INN (See the blue block in Fig.\ \ref{PI-INN}), it not only needs to learn the solutions of the forward problem but also needs to constrain the joint distribution $q(\bm{c},z)$. Here, $q(\bm{c},z)$ needs to satisfy two important assumptions \cite{INN}: $z$ must strictly follow the standard Gaussian distribution, and $\bm{c}$ and $z$ are independent of each other, i.e., $q(z|\bm{c}) = q(z)$. This independence ensures that redundant information is not encoded in both $\bm{c}$ and $z$. To achieve these assumptions, Ardizzone et al.\cite{INN} utilized MMD to minimize the distance between the joint distribution of output variables given by INN and the product of the marginal distribution given by data (see \ref{appendix-3} for details). While MMD can effectively distinguish between two distributions in low dimensions, it becomes less effective and requires a larger number of samples as the dimensionality increases \cite{Decreasing-Power-MMD}. In this study, a novel independence loss term is proposed, which dose not require labeled data for $\bm{c}$ and is applicable to high-dimensional cases, and the superiority of the loss is demonstrated in example \ref{experiment1}. The independence loss is defined as
\begin{equation}\label{latent-loss}
\mathcal{L}_{\mathrm{ind}} = \dfrac{1}{N}\sum\limits_{i=1}^{N}\big\Vert\mathrm{log}\,q(\hat{\bm{c}}^{(i)},\hat{z}^{(i)}) - \mathrm{log}\,q(\hat{\bm{c}}^{(i)},z^{(i)}) - (\mathrm{log}\,p(\hat{z}^{(i)}) - \mathrm{log}\,p(z^{(i)}))\big\Vert^2,
\end{equation}
Here $q(\bm{c},z)$ is the joint distribution given by INN, which can be obtained from the prior distribution $p(\lambda)$ of $\lambda$ and $g_\theta$, $[\hat{\bm{c}}^{(i)},\hat{z}^{(i)}] = g_{\bm{\theta}}(\lambda^{(i)})$ with $\lambda^{(i)}\sim p(\lambda)$, $z^{(i)}$ is independently drawn from the standard Gaussian prior $p(z)$. As demonstrated in \ref{appendix-1}, the loss is nonnegative, and is zero iff the two assumptions hold.

The expression of the complete loss function is given as follows:
\begin{equation}\label{loss-function}
\mathcal{L} = \alpha\mathcal{L}_{\mathrm{equ}} + \beta \mathcal{L}_{\mathrm{bound}} + \gamma \mathcal{L}_{\mathrm{ind}},
\end{equation}
where
\begin{align}
\mathcal{L}_{\mathrm{equ}} = \dfrac{1}{N}\sum\limits_{i=1}^{N}\big\Vert\mathcal{N}(\hat{u}^{(i)};\lambda^{(i)})\big\Vert^2,\quad
\mathcal{L}_{\mathrm{bound}} = \dfrac{1}{N}\sum\limits_{i=1}^{N}\big\Vert\mathcal{B}(\hat{u}^{(i)})\big\Vert^2.
\end{align}
Here, $\mathcal{L}_{\mathrm{equ}}$ represents equation loss term, $\mathcal{L}_{\mathrm{bound}}$ is boundary loss term, and $\mathcal{L}_{\mathrm{ind}}$ is independence loss term. $\hat{u}^{(i)} = \sum_{j=1}^{P}c_j(\lambda^{(i)})\phi_j$ is the predicted solutions of forward problem, and $\alpha$, $\beta$, $\gamma$ are relative weights of three loss terms. The detailed training process is summarized in Algorithm \ref{training-algorithm}.
\begin{algorithm}[h]
	\caption{Training procedure for PI-INN}
	\vspace{+2pt}
	\hspace*{0.001in} {\bf Input:} 
	Training set: $\{\lambda^{(i)}\}_{i=1}^{N}$, number of epochs: $E_{\mathrm{train}}$, learning rate: $\eta$, batch size: $N_{\mathrm{batch}}$, loss function: $\mathcal{L}$, weights of loss terms: $\alpha, \beta, \gamma$, coordinates set: $\{\bm{x}^{(j)}\}_{j=1}^{M}$, number of terms in (\ref{solution-expansion}): $N_P$.\\
	
	\vspace{-8pt}
	\begin{algorithmic}
		\For{epoch = $1:E_{\mathrm{train}}$}
		\State Sample minibatches from the training data and the latent space: $\{\lambda_{i}\}_{i=1}^{N_{\mathrm{batch}}}$, $\{z_i\}_{i=1}^{N_{\mathrm{batch}}}$, $z_i\stackrel{\mathrm{i.i.d}}{\sim}\mathcal{N}(0,\bm{I})$;
		\State Calculate $[\hat{\bm{c}}^{(i)}, \hat{z}^{(i)}] = g_{\bm{\theta}}(\lambda^{(i)})$;\Comment{$g_{\bm{\theta}}$ is the forward map of invertible neural network.}
		\State Substitute $\{\bm{x}^{(j)}\}_{j=1}^{M}$ into the NB-Net, and assemble matrix $\Phi\in\mathbb{R}^{M\times P}$;
		\State Let $\hat{u}^{(i)} = \Phi\hat{\bm{c}}^{(i)}$;
		\State Calculate $L = \mathcal{L}(\lambda^{(i)},\hat{u}^{(i)},\hat{z}^{(i)},z^{(i)})$;
		\State $\nabla\bm{\theta}\leftarrow\mathrm{Backpropagation}(L)$;
		\State $\bm{\theta}\leftarrow\bm{\theta} - \eta\nabla\bm{\theta}$;
		\EndFor
	\end{algorithmic}
	\hspace*{0.001in} {\bf Output:}
	Trained PI-INN.
	
	\vspace{+4pt}
	\label{training-algorithm}
\end{algorithm}

In the process of inversion, Bayesian inference can be accomplished by running the INN reversely, as shown by the green line in Fig. \ref{PI-INN}. Given the measurement $\tilde{u}$, the unknown input parameter can be repeatedly sampled with the latent variables $z$ shaped as a Gaussian distribution, and $\bm{c}$ keeps constant. Additionally, the desired posterior $p(\lambda\,|\,\tilde{u})$ is represented by the PI-INN, which transforms the known distribution $p(z)$ to the unknown input parameter space conditioned by different observations. We have theoretically proved that a well-trained PI-INN can efficiently sample, and explicitly provide the true posterior distribution without multiple likelihood evaluations required by classical Bayesian inference methods (see \ref{appendix-1}). 

To solve the optimization problem in Eq.\,(\ref{lse-problem}) and obtain $\bm{c}$, we approximate the marginal distribution $p(\bm{c})$ relying on the independence between $\bm{c}$ and $z$:
\begin{equation}\label{log-P-c}
\mathrm{log}\,p(\bm{c})\approx \mathrm{log}\,q(\bm{c},z) - \mathrm{log}\,p(z),
 \end{equation}
The gradient descent method is applied to solve problem (\ref{lse-problem}), and the detailed inference procedure is illustrated in Algorithm \ref{test-algorithm}.

\begin{algorithm}[ht]
	\caption{Inference procedure for PI-INN}
	
	\vspace{+4pt}
	\hspace*{0.001in} {\bf Input:}
	Trained PI-INN, locations of sensors: $\{\tilde{\bm{x}}^{(i)}\}_{i=1}^{O}$, observations: $\tilde{\bm{u}}$, number of inversion samples: $N_{\mathrm{samples}}$, regularization parameter: $\rho$.\\
	
	\vspace{-8pt}
	\begin{algorithmic}
		\State
		Substitute $\{\tilde{\bm{x}}^{(i)}\}_{i=1}^{O}$ into the NB-Net and assemble matrix $\tilde{\Phi}\in\mathbb{R}^{O\times N_P}$;
		\State Let $\mathcal{L}_{reg} = \Vert\tilde{\bm{u}} - \tilde{\Phi}\bm{c}\Vert^2 - \rho\mathrm{log}\,p(\bm{c})$;
            \State Minimize $\mathcal{L}_{reg}$ by gradient descent and obtain $\tilde{\bm{c}}$;
		\For{$k$ = $1:N_{\mathrm{samples}}$}
		\State Draw $z^{(k)}\sim\mathcal{N}(0,I)$;
		\State  Calculate $\tilde{\lambda}^{(k)} = g_{\bm{\theta}}^{-1}(\tilde{\bm{c}},z^{(k)})$;
		\EndFor 
	\end{algorithmic}
	\hspace*{0.001in} {\bf Output:}
	Samples of posterior distribution: $\{\tilde{\lambda}^{(k)}\}_{k=1}^{N_{\mathrm{samples}}}$.
	
	\vspace{+4pt}
	\label{test-algorithm}
\end{algorithm}

\begin{remark}
	Although PI-INN is designed for unsupervised learning, it can also be applied to labeled data, where a unified framework for both data-driven learning and physics-informed learning is introduced in \ref{appendix-2}.
\end{remark}
%=============================================================================================================

\section{Numerical experiments}\label{section:experiments}
In this section, several experiments are given to confirm the effectiveness of the PI-INN model. Firstly, the performance of our independence loss term is demonstrated and compared with the MMD loss term in \cite{INN} through solving an inverse kinematics problem. Then, the PI-INN model is applied to solve the inverse problems of 1-d and 2-d diffusion equations, and the inversion results are compared with MCMC. Finally, the practical applicability of the PI-INN model is also demonstrated in the seismic traveltime tomography problem. In examples involving inverse problems of PDEs, $\mathcal{L}_{\mathrm{equ}}$ can be defined as the variational formulation, thereby circumventing the need for computing second-order derivatives, especially when ReLU activation functions are employed (see \ref{appendix-4}). The accuracy of posterior samples can be evaluated by the $L_2$ relative errors defined as follows:
\begin{align}\label{l2-err}
\text{Relative}\ L_2\ \text{error} = \dfrac{\Vert \tilde{\lambda} - \lambda\Vert_2}{\Vert \tilde{\lambda}\Vert_2},
\end{align}
where $\lambda$ and $\tilde{\lambda}$ is the predicted and reference value, respectively.

\subsection{Inverse kinematics}\label{experiment1}
\begin{figure}[htp]
	\centering
	\subfigure[]{
		\includegraphics[scale=0.50]{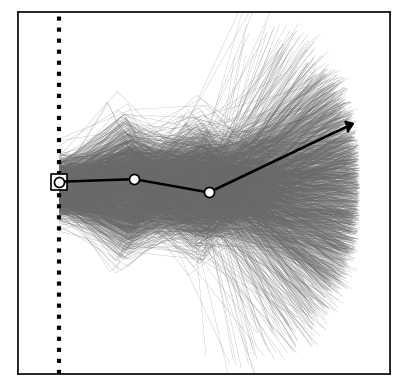}
		%\caption{fig1}
	}
	\subfigure[]{
		\includegraphics[scale=0.32]{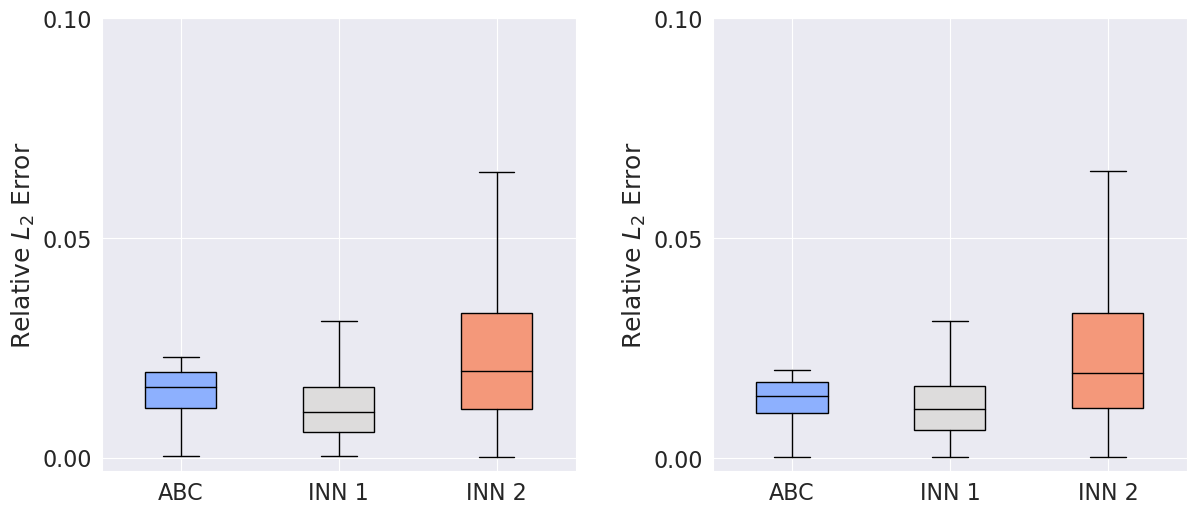}
	}
	
	\vspace{-5pt}
	\caption{Inverse kinematics: (a) Articulated arm model; (b) Boxplots of relative $L_2$ errors of reconstructing the position of the endpoints by three methods, and the left and right figures corresponding to Test 1 and Test 2, respectively.}
	\label{ex1-kinematics-boxplot}
\end{figure}

In this subsection, the accuracy and efficiency of our independence loss term will be demonstrated by the inverse kinematics problem introduced in \cite{INN}, where an articulated arm moves vertically along a rail and rotates at three joints. In Fig.\ \ref{ex1-kinematics-boxplot} (a), the dotted line symbolizes the rail, and the position coordinate of the arm is denoted by $x_1$, while the rotation angles of the three joints are represented by $x_2$, $x_3$, and $x_4$. The forward problem is to calculate the coordinate of end point $\bm{y}$ of the articulated arm. Corresponding map can be defined as follows:
\begin{equation}
\begin{array}{ll}
y_1 = l_1\,\mathrm{cos}(x_2) + l_2\,\mathrm{cos}(x_3 - x_2) + l_3\,\mathrm{cos}(x_4 - x_2 - x_3),\\
y_2 = x_1 +\ l_1\,\mathrm{sin}(x_2) + l_2\,\mathrm{sin}(x_3 - x_2)
+ l_3\,\mathrm{sin}(x_4 - x_2 - x_3).
\end{array}
\end{equation}
where $l_1 = 0.5$, $l_2 = 0.5$, and $l_3 = 1.0$ are the length of three arms. The prior of parameters are specified by normal distributions: $x_i\sim\mathcal{N}(0, \sigma_i)$, and $\sigma_1 = 0.25,\ \sigma_2 = \sigma_3 = \sigma_4 = 0.5\,\mathrm{rad} \ \hat{=}\ 28.65^{\circ}$. The inverse problem is to determine all the feasible parameters that generate the endpoint $\tilde{\bm{y}} = (\tilde{y}_1, \tilde{y}_2)$ of the arm.

To evaluate the efficiency of the independence loss term, we trained two data-driven INN models, referred to as INN 1 and INN 2, for solving the inverse problem (NB-Net was not employed in this particular experiment). Given training set $\{(\bm{x}^{(i)}, \bm{y}^{(i)})\}_{i=1}^{N}$, where $\bm{x} = (x_1,x_2,x_3,x_4)$, $\bm{y} = (y_1,y_2)$, the different loss functions are defined for two INN models as:
% MMD and our independence loss terms are leveraged to constrain the output distributions of INN 1 and INN 2, respectively. 
\noindent Loss functions of INN 1:
\begin{equation}\label{model1-loss}
\mathcal{L} = \dfrac{1}{N}\sum\limits_{i=1}^{N}\underbrace{\alpha\big\Vert \bm{y}^{(i)} - \hat{\bm{y}}^{(i)}\big\Vert^2}_{\text{\uppercase\expandafter{\romannumeral1}}} + \beta\underbrace{\big\Vert\mathrm{log}\,q(\hat{\bm{y}}^{(i)},\hat{z}^{(i)}) - \mathrm{log}\,q(\hat{\bm{y}}^{(i)},z^{(i)}) - (\mathrm{log}\,p(\hat{z}^{(i)}) - \mathrm{log}\,p(z^{(i)}))\big\Vert^2}_{\text{\uppercase\expandafter{\romannumeral2}}},
\end{equation}
\noindent Loss functions of INN 2:
\begin{equation}\label{model2-loss}
\begin{split}
\mathcal{L} = \dfrac{1}{N}\sum\limits_{i=1}^{N}\underbrace{\alpha\big\Vert \bm{y}^{(i)} - \hat{\bm{y}}^{(i)}\big\Vert^2}_{\text{\uppercase\expandafter{\romannumeral1}}} &+ \underbrace{\beta\mathrm{MMD}\big[p(\bm{y}^{(i)})p(z^{(i)}),q(\bm{y}^{(i)},z^{(i)})\big]}_{\text{\uppercase\expandafter{\romannumeral2}}} + \underbrace{\gamma\mathrm{MMD}\big[p(\bm{x}^{(i)}),q(\bm{x}^{(i)})\big]}_{\text{\uppercase\expandafter{\romannumeral3}}},
\end{split}
\end{equation}
where $\hat{\bm{y}}^{(i)}$ is the predicted value for $\bm{y}^{(i)}$, $\alpha, \beta, \gamma$ represent the weights of different loss terms, respectively. Part $\text{\uppercase\expandafter{\romannumeral1}}$ in Eqs. (\ref{model1-loss}) and (\ref{model1-loss}) is the forward-fitting loss term. In Eq. (\ref{model1-loss}), $\text{\uppercase\expandafter{\romannumeral2}}$ is the independence loss. In Eq. (\ref{model2-loss}), $\text{\uppercase\expandafter{\romannumeral2}}$ represents the MMD loss, and a reverse MMD loss $\text{\uppercase\expandafter{\romannumeral3}}$ is added to accelerate convergence, as in \cite{INN}. The inverse multiquadratic $k(x,x^{'}) = 1 / (1 + \Vert(x - x^{'})/h\Vert^2)$ is used as the kernel function of MMD, where $h = 1.2$.

\begin{figure}[htp]
	\centering
	\subfigure{
		\includegraphics[scale=0.265]{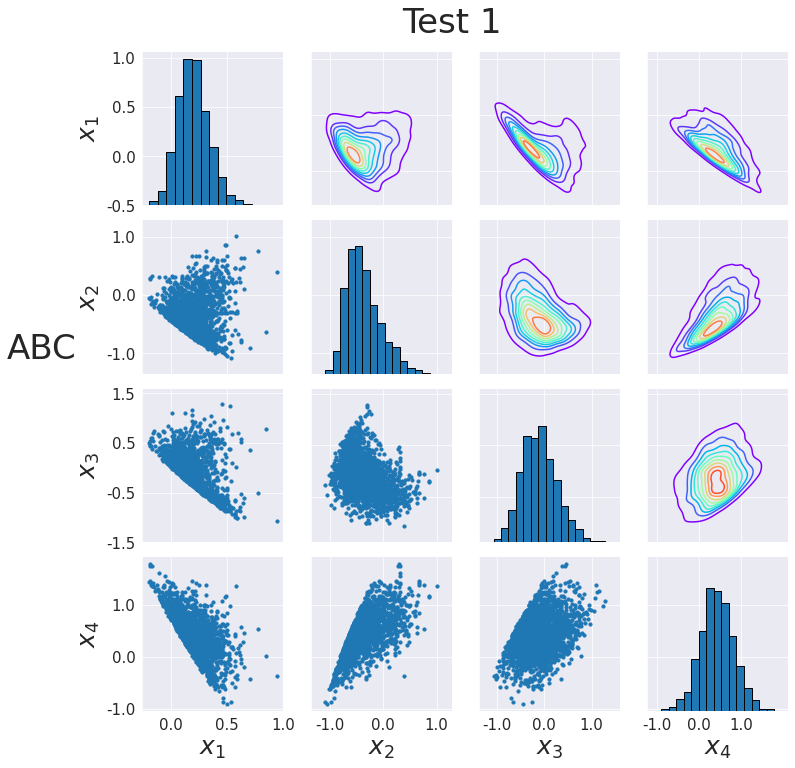}
		%\caption{fig1}
	}
	\subfigure{
		\includegraphics[scale=0.265]{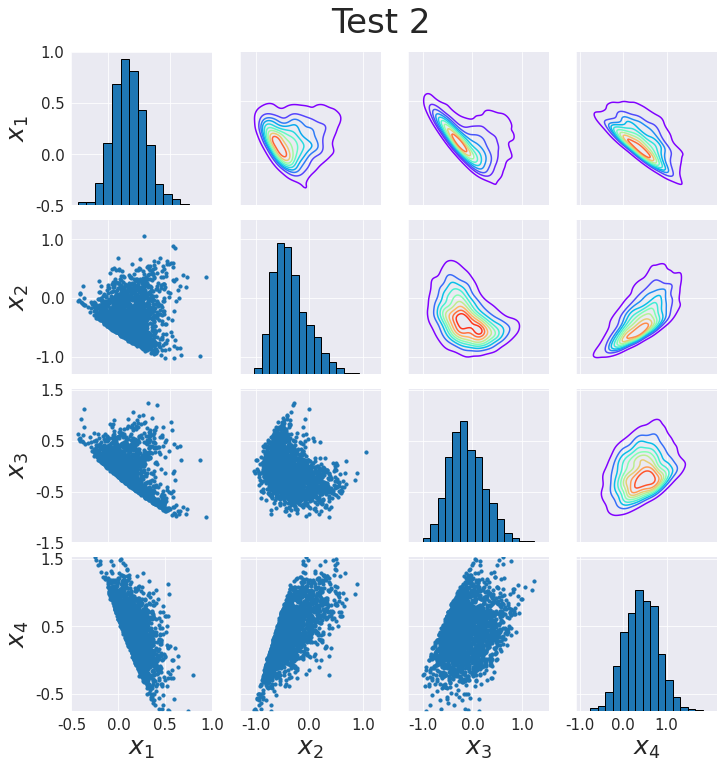}
	}\\
	\subfigure{
		\includegraphics[scale=0.265]{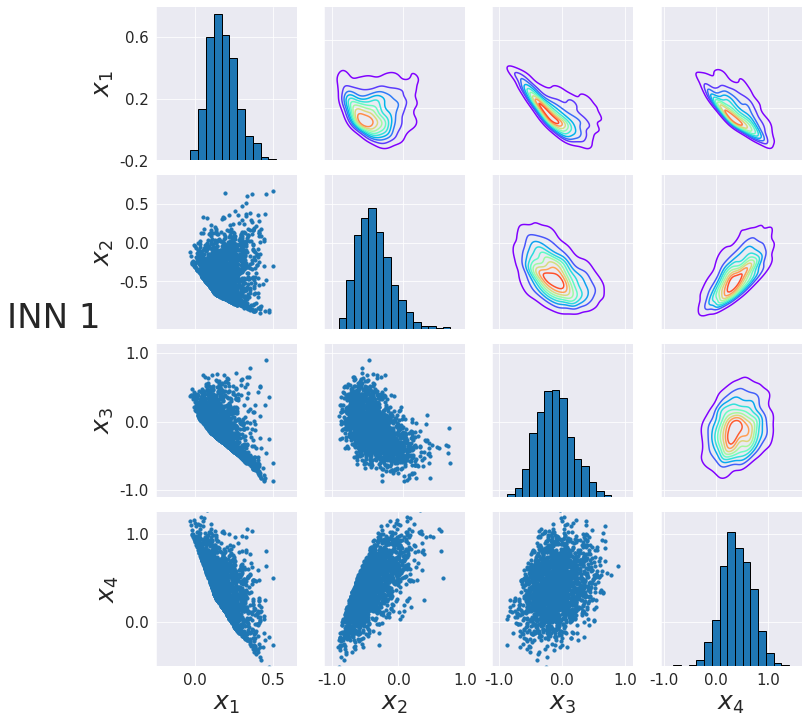}
		%\caption{fig1}
	}
	\subfigure{
		\includegraphics[scale=0.265]{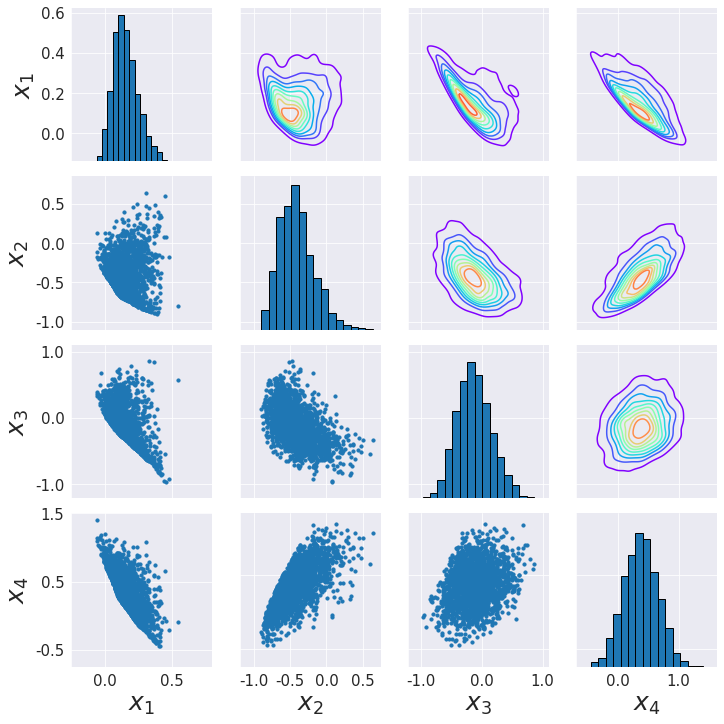}
	}\\
	\subfigure{
		\includegraphics[scale=0.265]{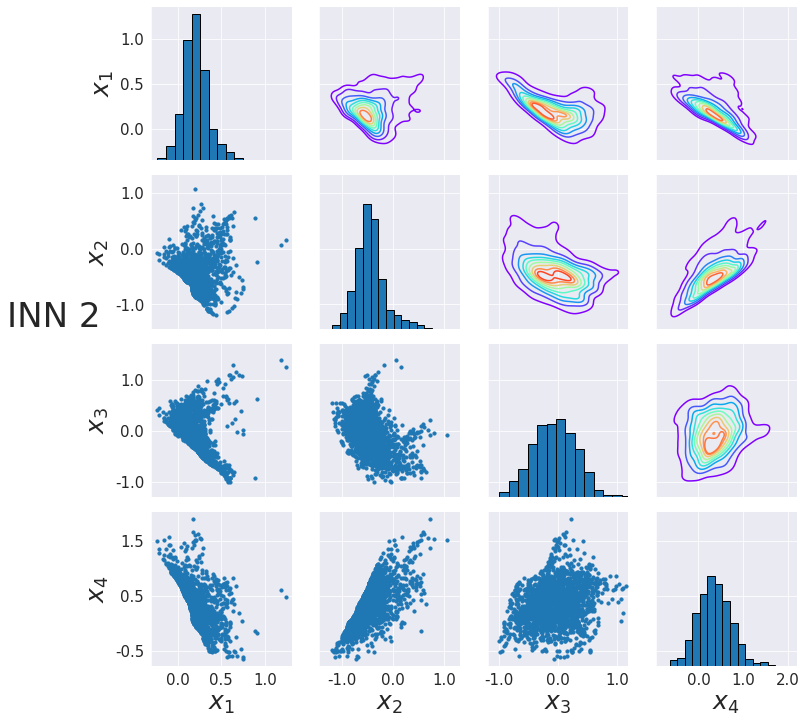}
		%\caption{fig1}
	}
	\subfigure{
		\includegraphics[scale=0.265]{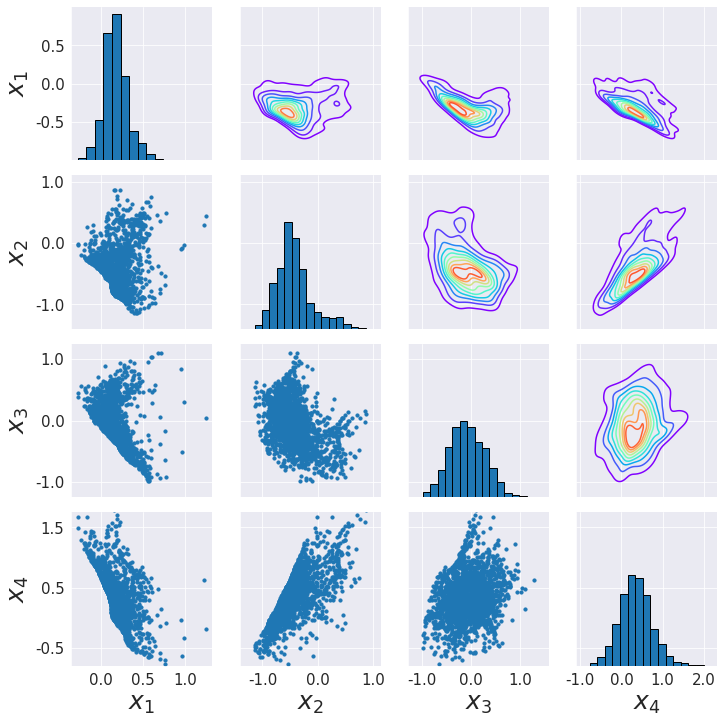}
	}\\
	
	\vspace{-10pt}
	\caption{Inverse kinematics: Sample distributions of two test cases for ABC and two INN models. In each figure, the diagonal display the histograms of four parameters, while the lower left triangle and upper right triangle areas show the pairwise scatter diagram and density estimation of the same four parameters, respectively.}
	\label{ex1-sample-distribution}
\end{figure}

Because the likelihood becomes computationally intractable as the test measurement of $\bm{y}$ is noise-free,  the approximate Bayesian computation (ABC) is utilized to obtain samples from the true posterior with comparison to two INN models (see more details in \ref{appendix-4-1}) instead of MCMC. Fig.\ \ref{ex1-kinematics-boxplot} (b) presents the boxplots of the relative $L_2$ errors for reconstructing the positions of the endpoints of two test cases. It can be seen that INN 1 outperforms INN 2, as it exhibits a lower mean and variance of relative errors. Moreover, the samples of posterior distributions for two test cases are displayed in Fig.\,\ref{ex1-sample-distribution}. It is can be observed that the sample distributions of INN 1 are closer to that of ABC than INN 2, and our independence loss term exhibits better performance than MMD loss.

% \textcolor{red}{We note that the INN model in \cite{INN} is trained with $600,000$ data, whereas in our example, a more efficient training strategy with only $4000$ data is employed. It's worth noting that the effectiveness of the MMD loss is influenced by the sample size \cite{Decreasing-Power-MMD}. Consequently, INN 2 may not be as efficient as the INN model described in \cite{INN}. Hence, we assert that our method exhibits superior performance in scenarios with smaller datasets.}

%%=============================================================================================================
\subsection{Application to 1-d diffusion equation for Gaussian and mixed Non-Gaussian random field}\label{experiment2}
We consider the following diffusion equation:
\begin{equation}\label{ex2-equ}
\left\{
\begin{array}{ll}
-\dfrac{\mathrm{d}}{\mathrm{d}x}\big(D(x)\dfrac{\mathrm{d}}{\mathrm{d}x}u(x)\big) = 5,\quad x\in [0, 1],\\
u(0) = 0,\quad u(1) = 1,
\end{array}\right.
\end{equation}
where $u(x)$ represents the concentration and $D(x)$ is the diffusion coefficient. In this example, we assume $D(x;\omega) = \mathrm{exp}(\tilde{D}(x;\omega))$ follows an exponential form of a Gaussian random field (GRF), where $\tilde{D}$ is defined in Eq. (\ref{ex2-gp-field}). The Karhunen-Lo$\grave{\text{e}}$ve expansion (KLE) of $\tilde{D}(x;\omega)$ is used to generate samples of $D(x;\omega)$, and the finite element method (FEM) is applied to obtain reference solutions with a uniform grid and the degree of freedom $N_{\text{Dof}} = 201$. The observation data is obtained by reference solutions with noise on sparse grid points.

\begin{equation}\label{ex2-gp-field}
\tilde{D}(x;\omega)\sim\mathcal{GP}\bigg(\frac{x}{2},\ \frac{9}{25}\text{exp}\Big(-\dfrac{\Vert x - x^{'}\Vert_2}{2}\Big)\bigg),\quad x,x^{'}\in [0, 1],
\end{equation}

To evaluate the effectiveness of our model, we conducted experiments using PI-INN to generate posterior samples of $D(x)$, where the importance sampling technique is utilized to enhance the precision of sampling (see \ref{appendix-IP}).
To better align with real-world scenarios, we introduced Gaussian noise with a standard deviation of $0.01$ to the measurements. For comparison, we employed MCMC as the reference method to obtain samples from the exact Bayesian posterior distribution.

\begin{figure}[htp]
	\centering
	\includegraphics[scale=0.325]{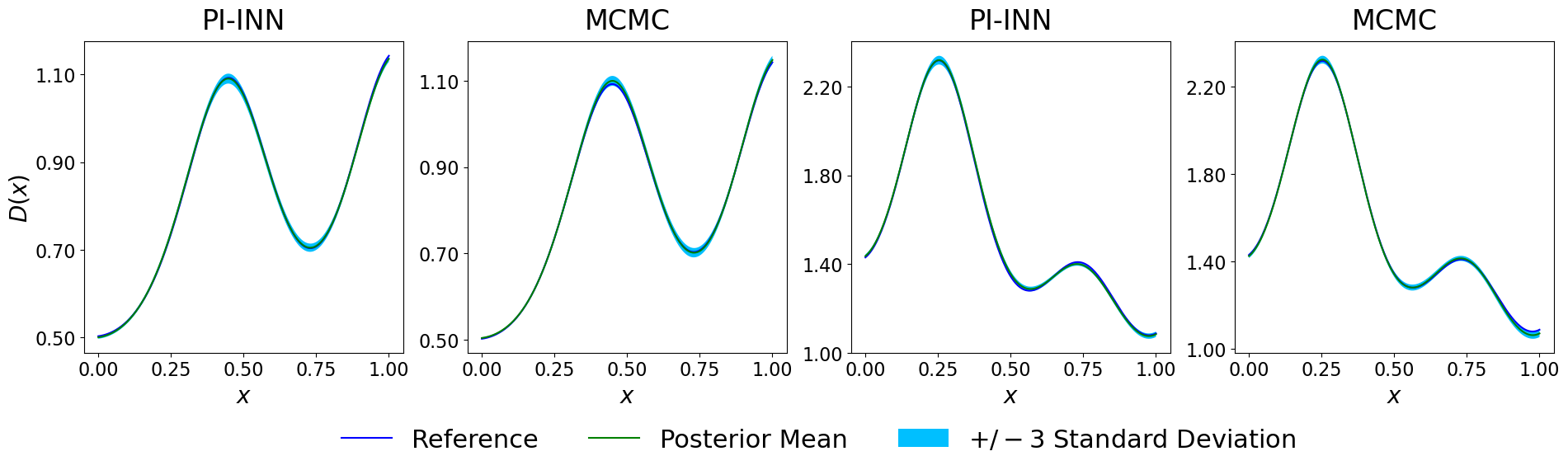}
	
	\vspace{-5pt}
	\caption{Two different test cases for 1-d diffusion equation under GRF prior: Comparisons of inversion results for PI-INN and MCMC, where the blue line and the green line represent the reference solution and the posterior mean, respectively. The shaded blue region represents the interval of the posterior mean $\pm$ 3 standard deviations, considering 11 evenly spaced sensors and Gaussian noise with a standard deviation of 0.01.}
	\label{ex2-diffusion-guassian-slices}
\end{figure}

Fig. \ref{ex2-diffusion-guassian-slices} illustrates the comparison of posterior samples obtained by PI-INN and MCMC across two test cases with different $u$. To quantity the effect of distribution approximation, the mean absolute point-wise errors for the mean and standard deviation of posterior distributions (see \ref{appendix-4-pointwise-error} for computation formula) are presented in Table~\ref{error-table}. Notably, in both test cases, the posterior samples generated by PI-INN closely align with the reference samples acquired through the MCMC approach. 

To further demonstrate the performance of our model, we consider the following mixed Non-Gaussian random field as the prior of $D(x)$:

\vspace{-5pt}
\begin{equation}\label{ex2-mixed-ngp-field}
D(x,\omega) = \mathrm{exp}\Big(\tilde{D}(x,\omega) + \frac{3}{4}m(x,\omega)\Big),
\end{equation}
where
\begin{equation*}\label{ex2-Gaussian}
\tilde{D}(x;\omega)\sim\mathcal{GP}\bigg(2x(1 - x),\ \dfrac{9}{25}\text{exp}\Big(-\dfrac{\Vert x - x^{'}\Vert_2}{2}\Big)\bigg),\quad x,x^{'}\in [0, 1],
\end{equation*}
\begin{equation*}\label{ex2-mixed}
m(x,\omega) =
\begin{cases}
\mathrm{sin}(\dfrac{\pi x}{2}), &\mbox{if $\xi(\omega) \geq 0$},\\
-\mathrm{sin}(\dfrac{\pi x}{2}), &\mbox{if $\xi(\omega) < 0$}.
\end{cases}
\end{equation*}

\begin{figure}[htp]
	\centering
	\includegraphics[scale=0.325]{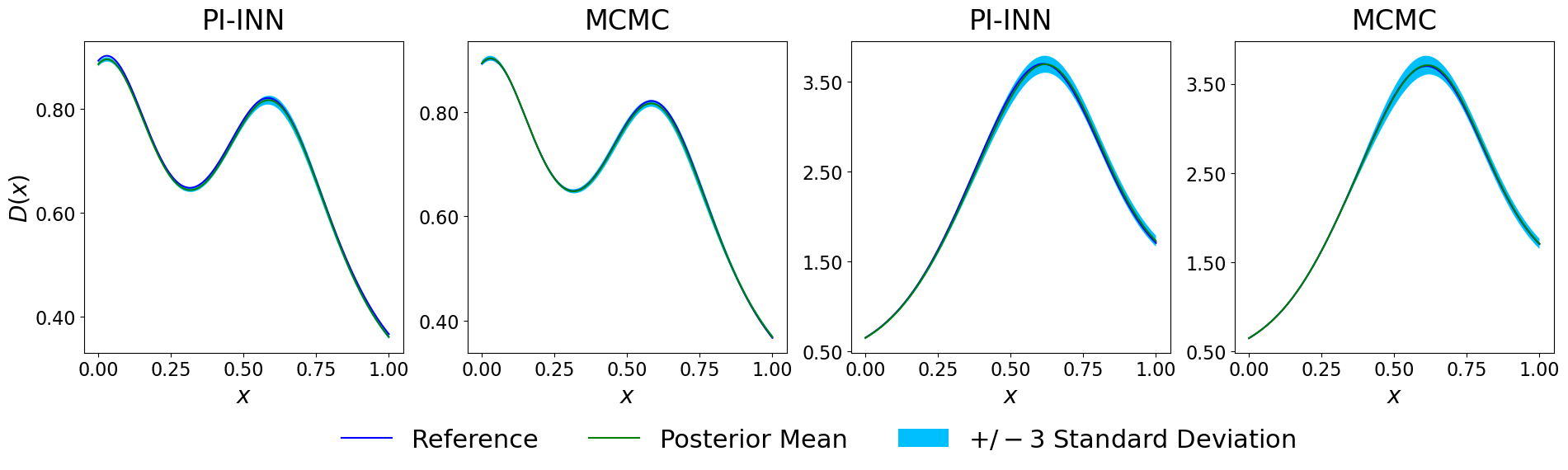}
	
	\vspace{-4pt}
	\caption{Two different test cases for 1-d diffusion equation under mixed Non-Gaussian prior: Comparisons of inversion results for PI-INN and MCMC, where the blue line and the green line represent the reference solution and the posterior mean, respectively. The shaded blue region represents the interval of the posterior mean $\pm$ 3 standard deviations, considering 11 evenly spaced sensors and Gaussian noise with a standard deviation of $0.01$.}
	\label{ex2-diffusion-mixed-guassian-slices}
\end{figure}

 The inversion results of two test cases are displayed the in Fig. \ref{ex2-diffusion-mixed-guassian-slices}. We also calculate the mean point-wise errors in Table\,\ref{error-table}. The results reveal a strong agreement between the posterior mean and standard deviation produced by our method and the reference solutions from MCMC, underscoring the efficiency and reliability of our approach.

\begin{table}[htbp]
	\centering
	\caption{Error table for posterior sample by PI-INN to approximate the reference distribution by MCMC (corresponding to the posterior samples in Fig.\,\ref{ex2-diffusion-guassian-slices} and Fig.\,\ref{ex2-diffusion-mixed-guassian-slices}).}
	\begin{tabular}{|c|ll|ll|}
		\hline
		\multirow{2}{*}{Prior} & \multicolumn{2}{c|}{Test 1}     & \multicolumn{2}{c|}{Test 2}     \\ \cline{2-5} 
		& \multicolumn{1}{c|}{Mean} & \multicolumn{1}{c|}{Std} & \multicolumn{1}{c|}{Mean} & \multicolumn{1}{c|}{Std} \\ \hline
		GRF                    & \multicolumn{1}{l|}{$4.26\times 10^{-3}$}     &  $3.64\times 10^{-4}$   & \multicolumn{1}{l|}{$8.47\times 10^{-3}$}     &  $1.50\times 10^{-3}$   \\ \hline
		Mixed Non-Gaussian     & \multicolumn{1}{l|}{$8.13\times 10^{-3}$}     &  $3.03\times 10^{-3}$   & \multicolumn{1}{l|}{$4.57\times 10^{-3}$}     & $1.08\times 10^{-3}$    \\ \hline
	\end{tabular}
	\label{error-table}
\end{table}

\begin{figure}[htp]
	\centering
	\subfigure{
		\includegraphics[scale=0.35]{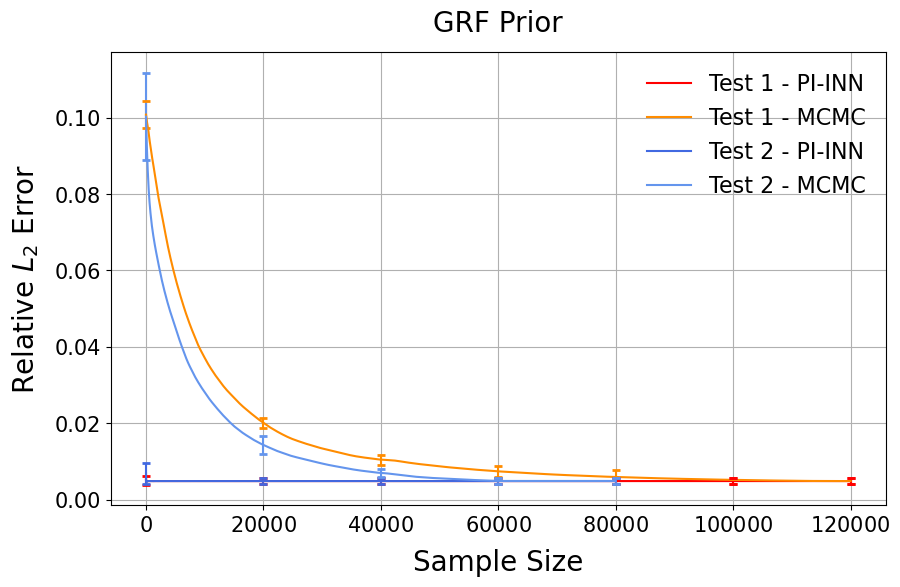}
	}
	\subfigure{
		\includegraphics[scale=0.35]{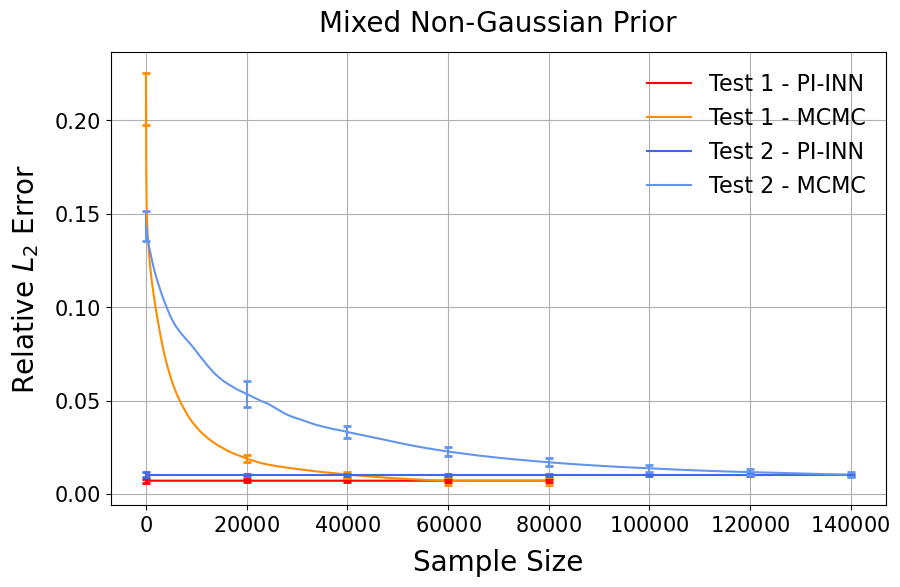}
	}
	
	\vspace{-10pt}
	\caption{Comparison of PI-INN and MCMC convergence using sample size vs. relative $L_2$ error for sample mean. Error bars represent results from 5 independent experiments. Total sample size for both subfigures: $80000$ (test 1) and $120000$ (test 2).}
	\label{ex2-efficiency}
\end{figure}

\begin{figure}[htp]
	\centering
	\subfigure{
		\includegraphics[scale=0.35]{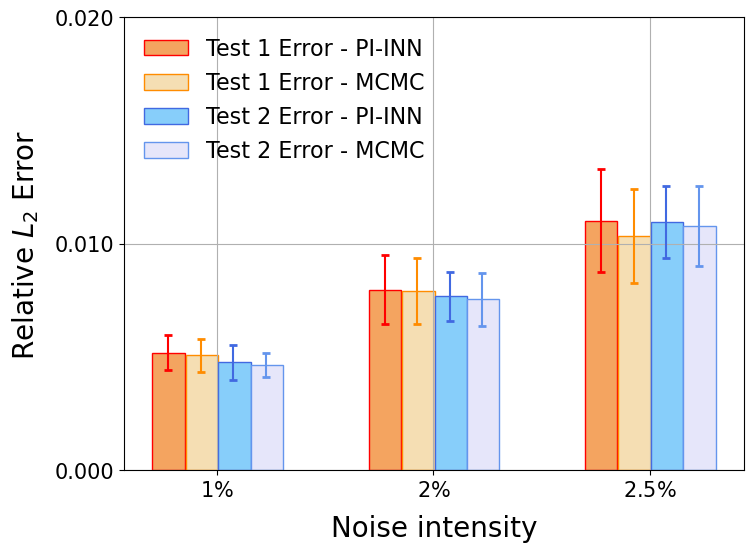}
		%\caption{fig1}
	}
	\subfigure{
		\includegraphics[scale=0.35]{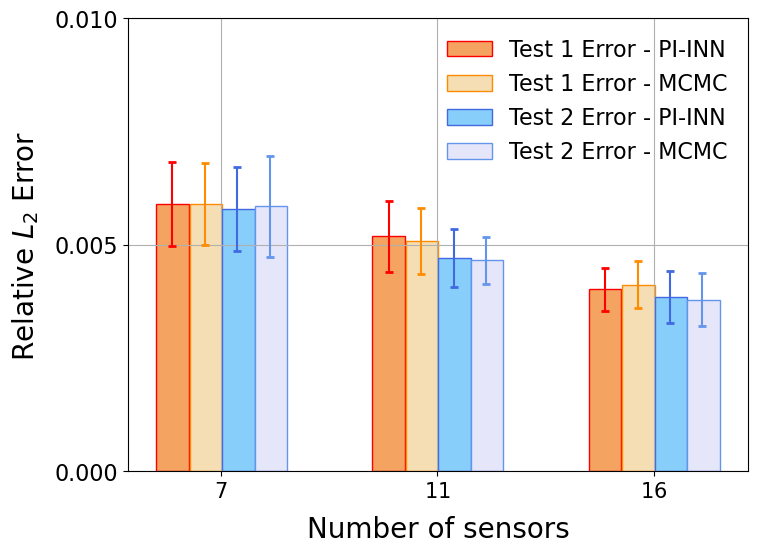}
	}
	\subfigure{
		\includegraphics[scale=0.35]{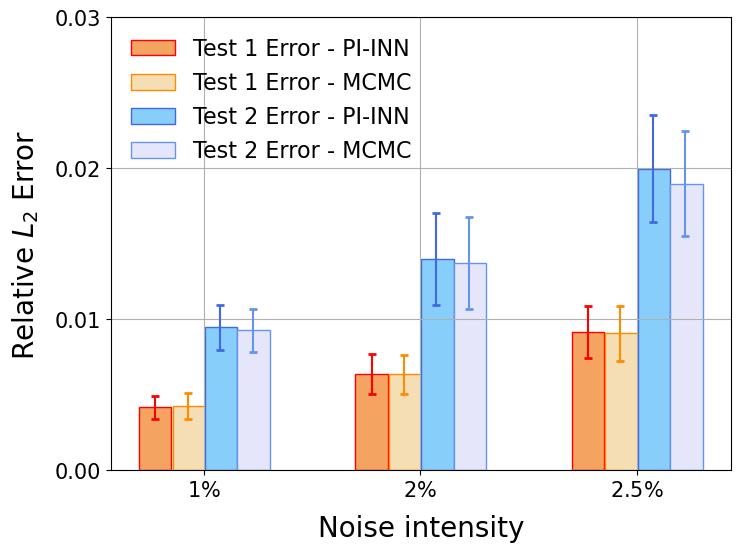}
		%\caption{fig1}
	}
	\subfigure{
		\includegraphics[scale=0.35]{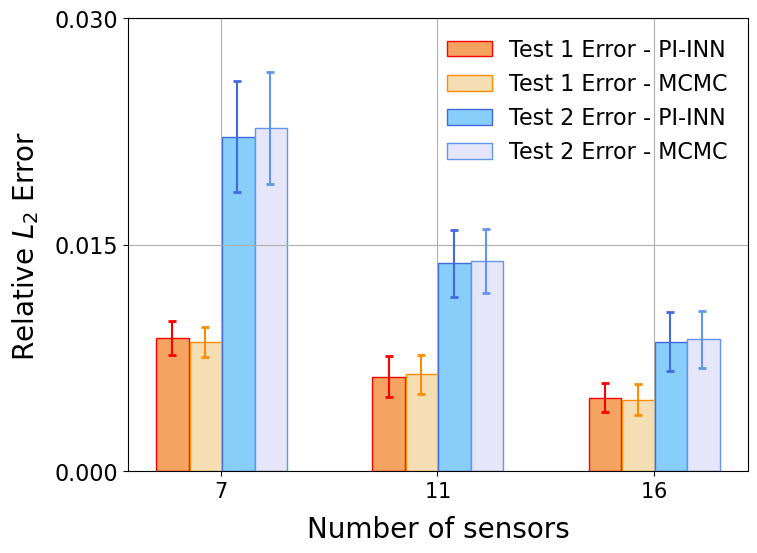}
	}
	
	\vspace{-10pt}
	\caption{Relative $L_2$ errors of $1000$ posterior samples with the standard deviations of noise and the number of sensors changed for two test cases of 1-d diffusion equation: Bars represent the averages and error bars show the standard deviations. The subfigures of the first and second row correspond to GRF and mixed Non-Gaussian random field prior cases, respectively.
 }
	\label{ex2-noise-sensor-analysis}
\end{figure}

Furthermore, the efficiency and robustness of our model is validated. For the two distinct prior assumptions considered in the aforementioned examples, Fig. \ref{ex2-efficiency} visually presents the convergence patterns of posterior samples obtained through PI-INN and MCMC. Remarkably, the posterior distribution samples derived via PI-INN exhibit a notably accelerated convergence rate when compared to those produced by MCMC. This disparity underscores the efficiency and computational advantages offered by the PI-INN approach. The robustness is validated by analyzing its performance under different combinations of measurement sensors and noise intensity. As shown in Fig.\,\ref{ex2-noise-sensor-analysis}, we observed that the mean and standard deviation of the relative errors of the posterior samples of PI-INN has similar order of magnitudes to those of MCMC under the different observation conditions. This indicates that the performance of our model is consistent with the MCMC method, which shows the robustness and reliability of PI-INN model.

%%=============================================================================================================
\subsection{Application to Darcy equation}\label{experiment3}
In this example, the Darcy equation in the two-dimensional domain is considered:
\begin{equation}\label{ex3-equ}
\left\{
\begin{array}{lll}
-\nabla\cdot(K(x)\nabla u(x)) = 0,\quad x\in[0,1]^2,\\
u(0, x_2) = 1,\ u(1, x_2) = 0,\\
\dfrac{\partial u}{\partial n}\Big|_{x_2 = 0} = \dfrac{\partial u}{\partial n}\Big|_{y = 1} = 0,
\end{array}\right.
\end{equation} 
where $u(x)$ represents the pressure field $K(x)$ is the permeability field, and can be sampled from $K(x,\omega) = \mathrm{exp}\big(\tilde{K}(x,\omega)\big)$, where $\tilde{K}(x,\omega)$ is defined by Eq. (\ref{ex3-K}). The domain can be divided into $64\times 64$ uniform mesh and the reference solutions can be obtained by FEM. $22\times 22$ evenly spaced grid points in $[0,1]^2$ are selected as the positions of sensors during the inversion process. The first $15$ leading terms of KLE are utilized as the input of INN.
\begin{equation}\label{ex3-K}
\tilde{K}(x,\omega) = \mathcal{GP}\Big(0,\ \text{exp}\big(-5\Vert x - x^{'}\Vert_2\big)\Big),\quad x,x^{'}\in [0,1]^2,
\end{equation}

\begin{figure}[ht]
	\centering
	\subfigure{
		\includegraphics[scale=0.34]{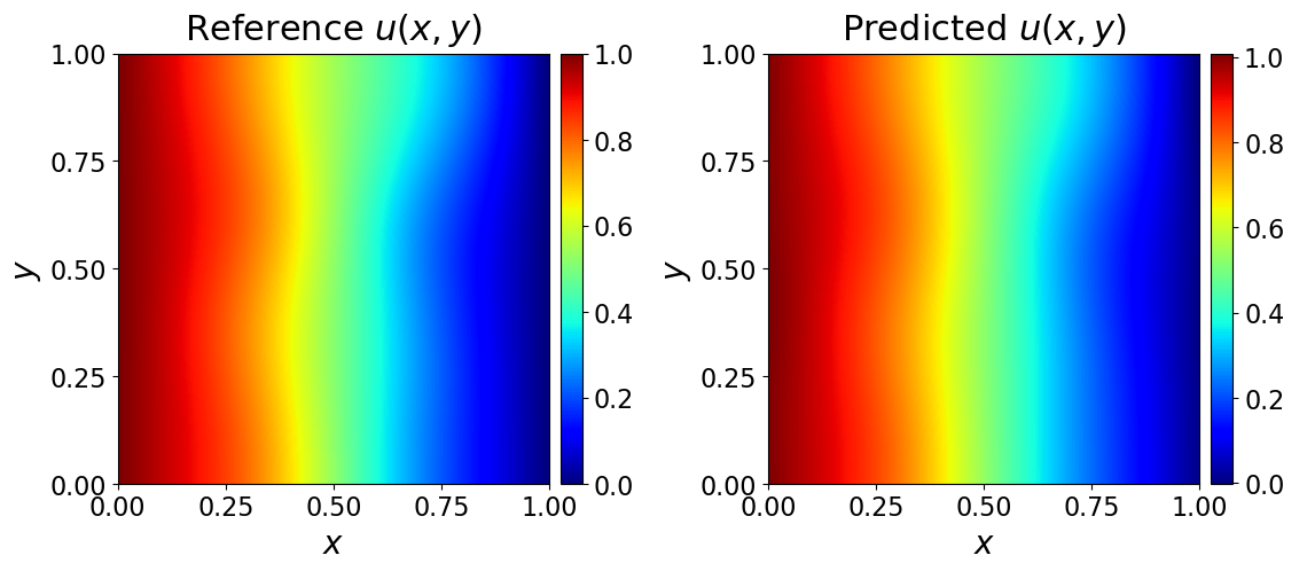}
	}\vspace{-10pt}
	\subfigure{
		\includegraphics[scale=0.34]{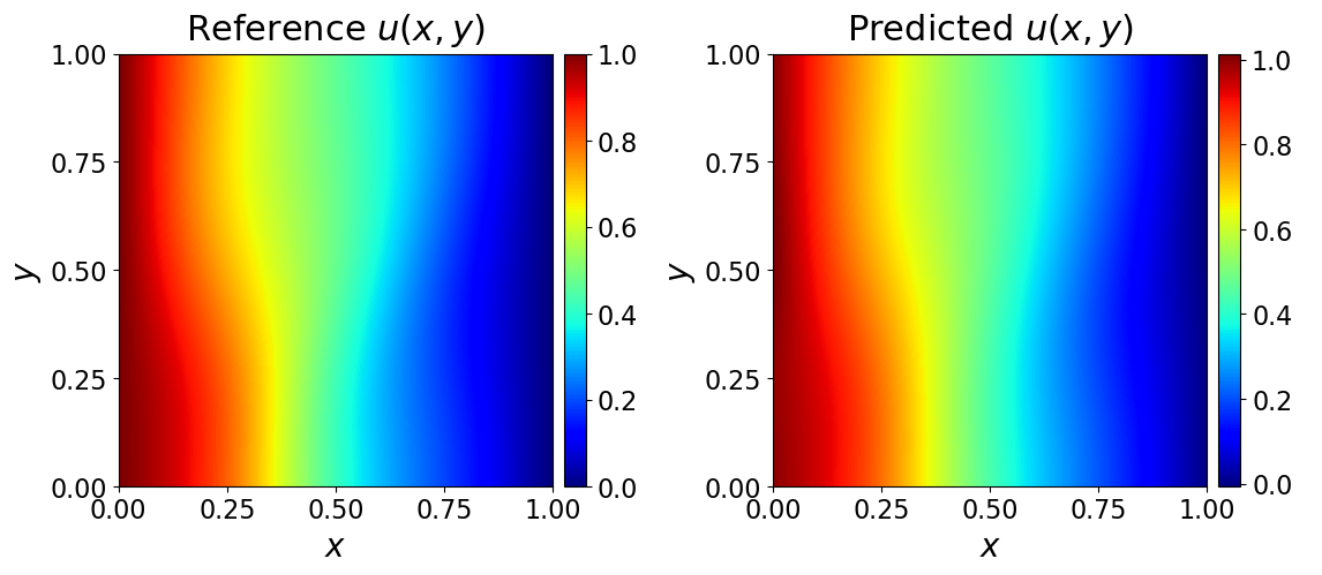}
	}
	\caption{Comparison of the reference solutions $u$ and the predicted solutions for 2-d Darcy equation by PI-INN model.}
	\label{ex3-darcy-simulations}
\end{figure}

\begin{figure}[htp]
	\centering
	\includegraphics[scale=0.32]{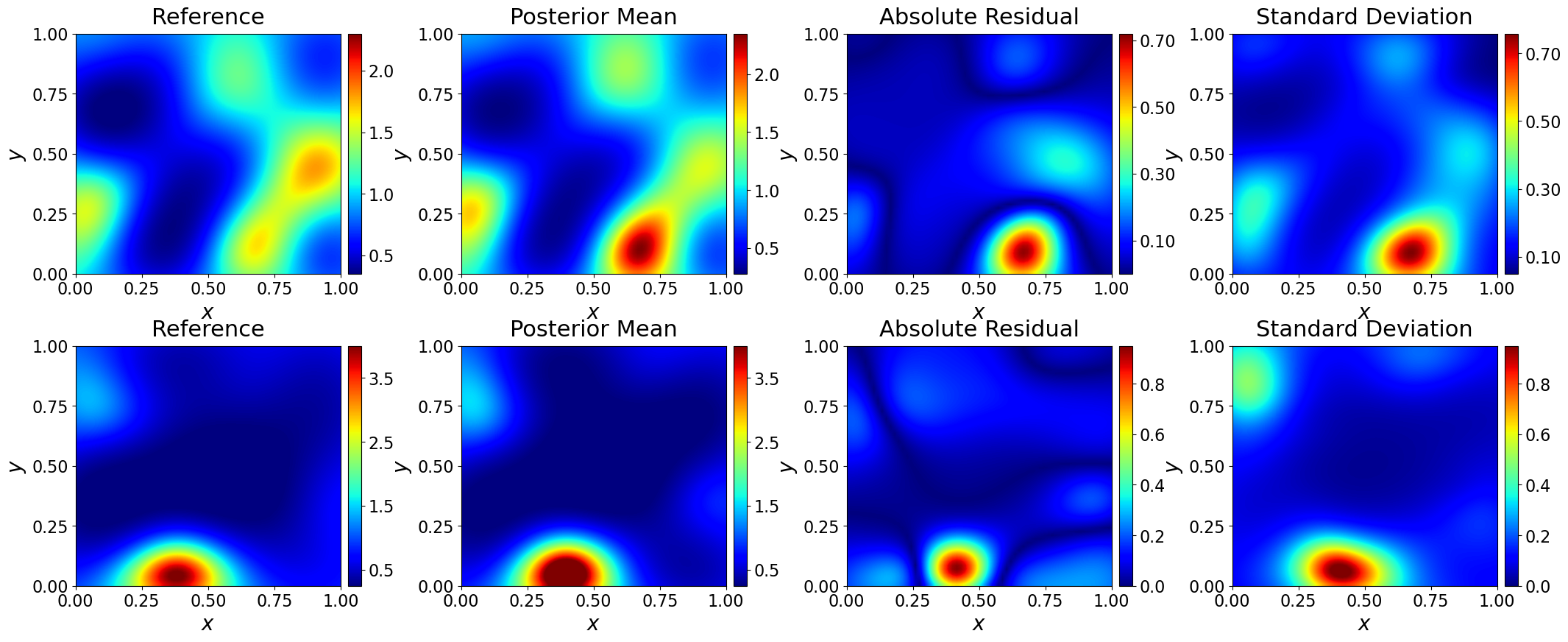}
	\caption{Results for Bayesian inference by PI-INN on two test cases: reference permeability fields (first column), mean and standard deviation of posterior samples (second and fourth columns), absolute errors between reference fields and posterior mean (third column). Zero-mean Gaussian noise with $0.01$ standard deviation and $22\times 22$ sensors are used in this example.}
\label{ex3-darcy-inversions}
\end{figure}

To assess the accuracy of our forward simulation, we display the predicted pressure fields by PI-INN for two test cases in Fig.\ \ref{ex3-darcy-simulations}. These solutions are almost identical to the reference solutions obtained by FEM, which indicates the high accuracy of our forward predictions. Fig.\,\ref{ex3-darcy-inversions} portrays the comparison between the posterior mean and reference permeability fields in two distinct test cases. We observe that the posterior samples generated by PI-INN can accurately capture the global features of the reference fields, while the standard deviation effectively quantifies the uncertainty.

\begin{figure}[htp]
	\centering
	\subfigure{
		\includegraphics[scale=0.37]{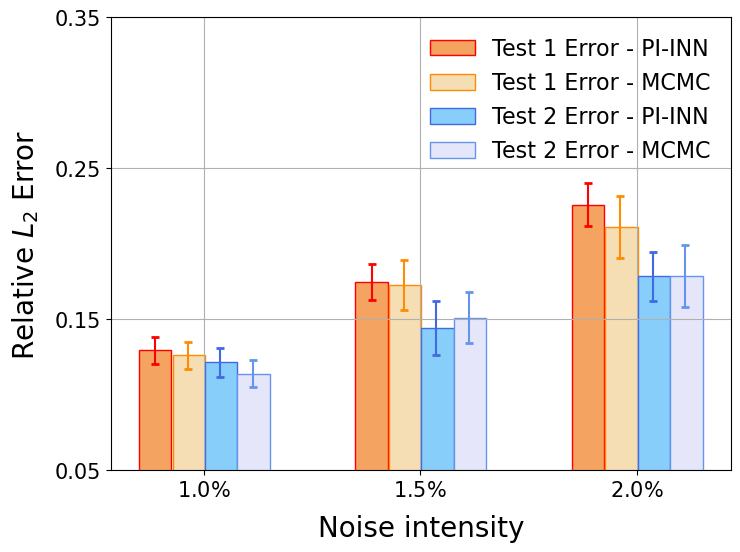}
		%\caption{fig1}
	}
	\subfigure{
		\includegraphics[scale=0.37]{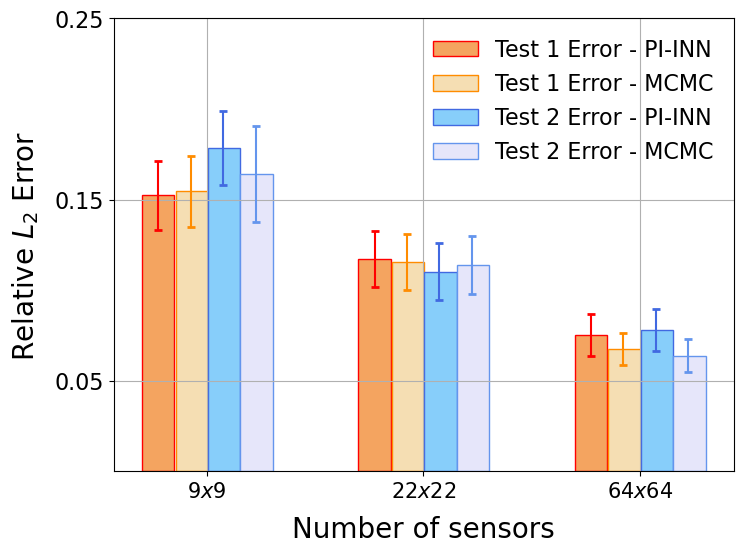}
	}
	\caption{Relative $L_2$ errors of $1000$ posterior samples with noise intensity and the number of sensors changed for 2-d darcy equation: Bars represent the averages and error bars show the standard deviations. $22\times 22$ sensors and zero-mean Gaussian noise with standard deviation $0.005$ are used in two subfigures, respectively.}
	\label{ex3-noise-sensor-analysis}
\end{figure}

To quantity the uncertainty of inversion results with the various observation conditions, we present the mean and standard deviation of the relative $L_2$ error of 1000 posterior samples in Fig.\ \ref{ex3-noise-sensor-analysis}. It is noted that the mean and standard deviation of the relative errors of PI-INN is consistent with those of MCMC, which demonstrates the effectiveness of our method.

%%=============================================================================================================
\subsection{Seismic traveltime tomography}\label{experiment4}
Seismic travel time tomography is a widely used technique to image the subsurface properties and interior structures of the Earth \cite{Tomo-NFs}. In seismic imaging, it is typically formulated as an inverse problem to estimate underground velocity models with the observed travel time data \cite{PINNtomo}. In this example, PI-INN is applied to solve the practical seismic tomography problem based on the first-arrival travel time.

Consider the 2-d Eikonal equation defined as follows:
\begin{equation}\label{eikonal-equ}
\left\{
\begin{array}{ll}
\vert\nabla T(\bm{x})\vert^2 = \dfrac{1}{v^2(\bm{x})},\ \bm{x}\in \mathcal{D},\\
T(\bm{x}_s) = 0,
\end{array}
\right.
\end{equation}
where $\mathcal{D} = [0, 4]^2$ represents the domain, and $T(\bm{x})$ is the travel time to any position $\bm{x}$ from the point-source $\bm{x}_s = (2, 0)$. $v(\bm{x})$ is the velocity model defined on $\mathcal{D}$, and $\nabla$ denotes the gradient operator. In this example, we assume that the velocity model can be parameterized as follows:
\begin{equation}\label{velocity-formula}
v(\bm{x}) = v(\bm{x}_0) + g_{X}(x - x_0) + g_{Y}(y - y_0),
\end{equation}
where
\begin{equation}\label{velocity-priors}
g_{X} = 0,\quad g_{Y} = \begin{cases}
g_{Y}^{1}, &\mbox{if\ $0 < Y \leq h_{1}$},\\
g_{Y}^{2}, &\mbox{if\ $h_{1} < Y \leq \sum\limits_{i=1}^{2}h_{i}$},\\
g_{Y}^{3}, &\mbox{if\ $\sum\limits_{i=1}^{2}h_{i} < Y \leq \sum\limits_{i=1}^{3}h_{i}$},\\
g_{Y}^{4}, &\mbox{if\ $\sum\limits_{i=1}^{3}h_{i} < Y \leq \sum\limits_{i=1}^{4}h_{i}$}.
\end{cases}
\end{equation}
Here, $v(\bm{x}_0)$ is the velocity on the origin $\bm{x}_0 = (0, 0)$, $g_{X}$ and $g_{Y}$ are the velocity gradients along the horizontal and vertical directions, respectively. In this setting, the variation of vertical velocity gradient is determined by the depth parameters $h_i (i=1, 2, 3, 4)$. Training data is generated with samples of the prior distributions of random parameters (such as vertical velocity gradient values and depths) in Eq. (\ref{params-prior}) and velocity models in Eq. (\ref{velocity-formula}). $6\times 6$ evenly spaced sensors are distributed to record the travel time data, which is calculated using fast sweeping method \cite{FSM} with $101\times 101$ uniform mesh. The source point is fixed at $(2, 0)$.

\begin{equation}\label{params-prior}
\begin{array}{cc}
g_{Y}^{1}\sim\mathcal{N}(0.2, 0.25),\ \ g_{Y}^{2}\sim\mathcal{N}(0.4, 0.25),\ g_{Y}^{3}\sim\mathcal{N}(0.5, 0.25),\ \ g_{Y}^{4}\sim\mathcal{N}(1, 1),\vspace{0.5ex}\\
h_{1}\sim\mathcal{U}[0.75, 1.25],\ \ h_{2}\sim\mathcal{U}[0.75, 1.25],\ \ 
h_{3}\sim\mathcal{U}[0.75, 1.25],\ \ h_{4} = 4 - h_{1} - h_{2} - h_{3},
\end{array}
\end{equation}

\begin{figure}[ht]
	\centering
	\subfigure{
		\includegraphics[scale=0.32]{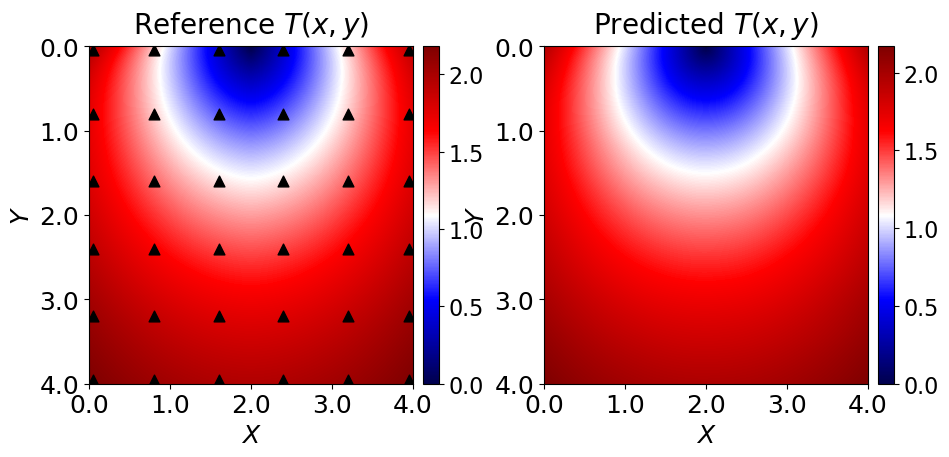}
	}\vspace{-10pt}
	\subfigure{
		\includegraphics[scale=0.32]{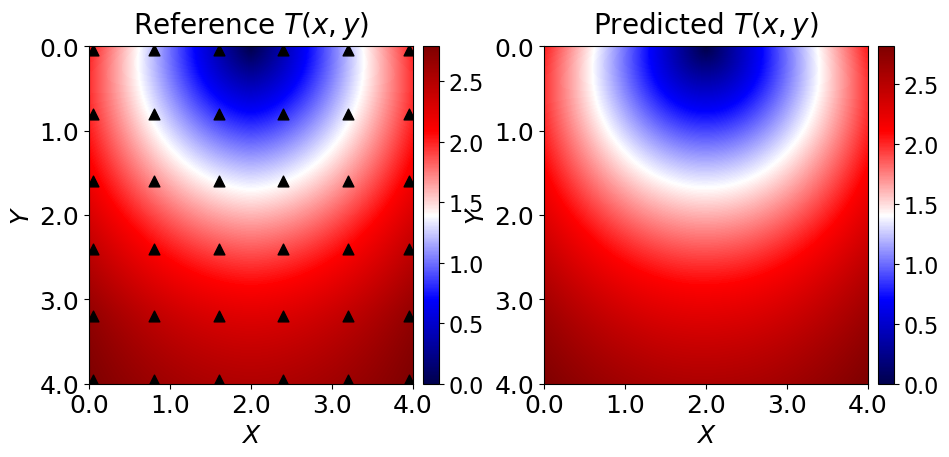}
	}
	\caption{Comparison of the reference travel time and the predicted ones by PI-INN for two test cases of seismic traveltime tomography. The triangle marks in the reference $T(x,y)$ represent the sensor locations for tomography.}
	\label{ex4-T-pred}
\end{figure}

\begin{figure}[ht]
	\centering
	\subfigure{
		\includegraphics[scale=0.29]{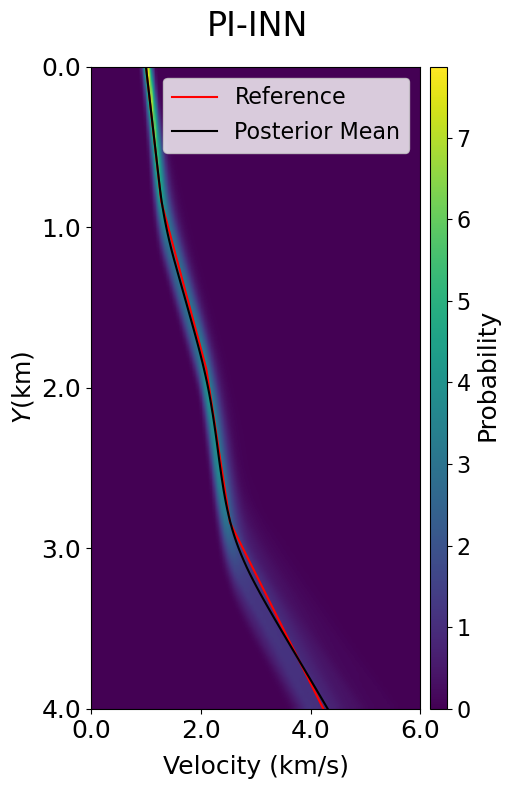}
	}
	\subfigure{
		\includegraphics[scale=0.29]{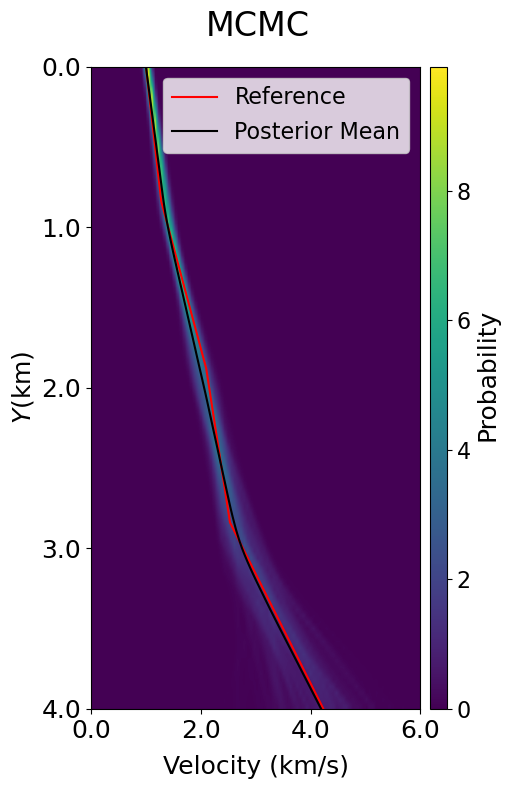}
	}
	\subfigure{
		\includegraphics[scale=0.29]{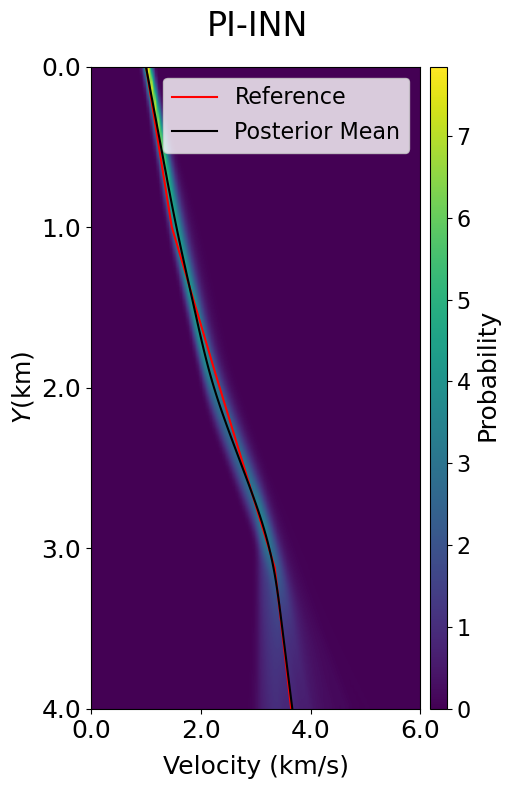}
	}
	\subfigure{
		\includegraphics[scale=0.29]{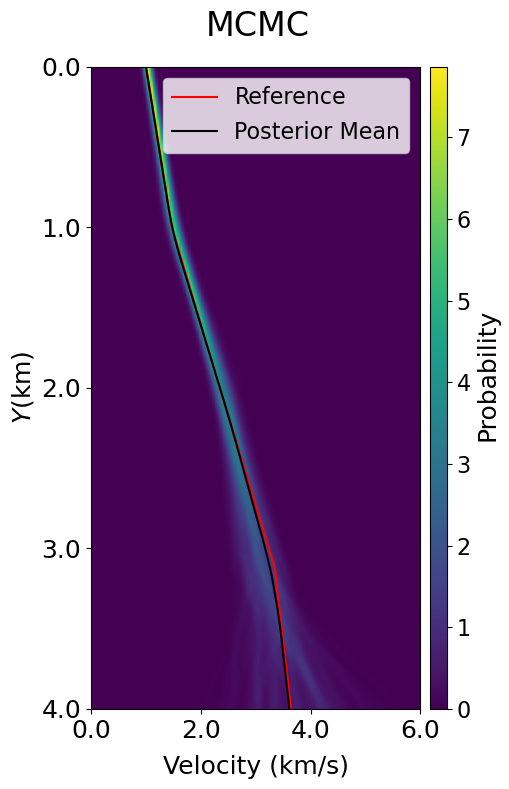}
	}
	\vspace{-10pt}
	\caption{The marginal posterior distributions obtained using PI-INN (first and third subfigures) and MCMC (second and fourth subfigures) for two test cases. $6\times 6$ sensors and zero-mean Gaussian noise with $0.1$ standard deviation are used in this example. Black and red lines show the reference solution and posterior mean velocity, respectively. The background color represents the kernel density estimation of posterior samples.}
	\label{ex4-INN-MCMC-tomo}
\end{figure}

\begin{figure}[htp]
	\centering
	\includegraphics[scale=0.46]{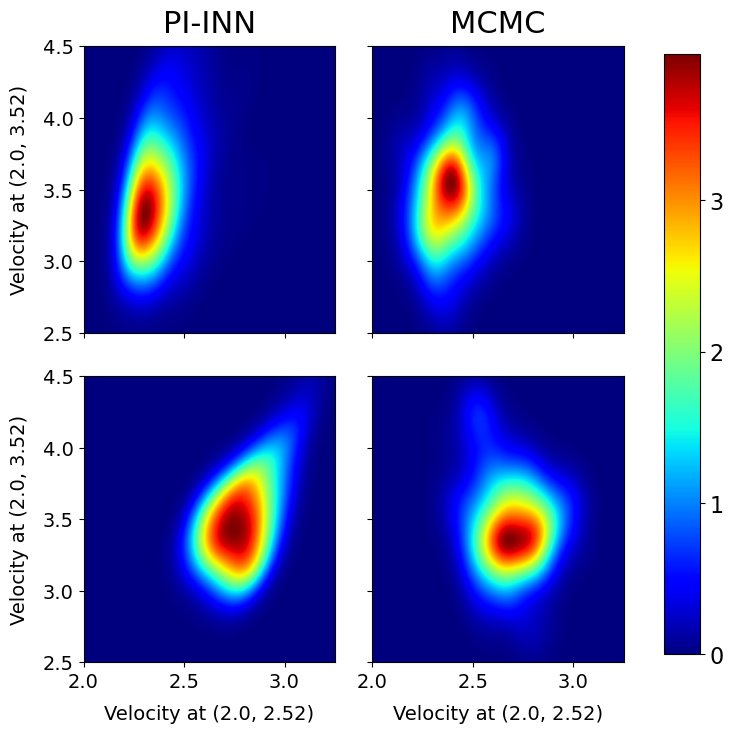}
	\caption{The bivariate marginal distributions of velocities at the locations of $(2.0, 2.52)$ and $(3.0, 3.52)$ obtained using PI-INN and MCMC, respectively, for the two velocity profiles in Fig.\ \ref{ex4-INN-MCMC-tomo}.}
	\label{ex4-bivariate-dist}
\end{figure}

As depicted in Fig.\ \ref{ex4-T-pred}, the predicted travel time fields by PI-INN are almost identical to the reference solutions, which indicates that 
PI-INN can provide high-accuracy forward predictions for Eq.\ (\ref{eikonal-equ}). To visualize the effect of travel time tomography, the comparison of marginal posterior distributions of PI-INN and MCMC along the vertical dimension is depicted in Fig.\ \ref{ex4-INN-MCMC-tomo}. Additionally, the bivariate marginal distributions of velocities for PI-INN and MCMC at two locations $(2.0, 2.52)$ and $(2.0, 3.52)$ are given in Fig.\ \ref{ex4-bivariate-dist} for the two test cases. The results demonstrate that PI-INN can efficiently recover the true velocity, quantifies uncertainty, and captures correlations of velocity models. These findings highlight the potential of our PI-INN to solve practical problems in geophysical inversion.

\section{Conclusions}\label{section:conclusions}
In this paper, a novel PI-INN is proposed to efficiently solve Bayesian inverse problems, which provides a general learning framework for Bayesian inverse problems. In contrast to traditional Bayesian inference, PI-INN can provide a tractable estimation of the posterior distribution, which enables efficient sampling and high-accuracy density evaluation through the novel independence loss term incorporated with physics-informed loss terms and is particularly advantageous for high-dimensional problems. Moreover, the asymptotic analysis of PI-INN is theoretically derived under the assumption that the loss function approaches zero. Numerical experiments show that the proposed PI-INN is stable and effective, and can also provide enough numerical accuracy for the indirect sparse and noisy measurement data, which supports the theoretical results of this paper. 

We will further improve the method and overcome the following limitations in future. First of all, it necessitates that $\mathrm{ndim}(\bm{c}) < \mathrm{ndim}(\lambda)$ within our framework. This requirement may not be suitable for scenarios where the dimensionality of $\lambda$ is low yet the solution remains complex.
%This constraint can be mitigated by introducing stochastic noise $\mathrm{v}$ at the input layer, i.e. the forward map is also probabilistic. In this modified setting, the dimensional relationship becomes $\mathrm{ndim}(\lambda) + \mathrm{ndim}(v) = \mathrm{ndim}(\bm{c}) + \mathrm{ndim}(z)$. But there is a lack of corresponding work to do that in the current framework.
Second, it's essential to note that our numerical experiments primarily incorporate observation noise characterized by low variance (see \ref{appendix-4-2-2} for details). In practical applications, it becomes imperative to address noise with higher degrees of uncertainty.
%Furthermore, while this work centers on inverse problems, a unified uncertainty qualification (UQ) framework could be extended. Current research available for reference includes \cite{PI-GANs, NFFs, PI-VAE, L-HYDRA, VI-NFs}.
Future work will focus on addressing the aforementioned limitations and exploring the various practical applications to solve Bayesian inverse problems in geophysics and Earth science.                                        
\section*{Data availability statement}
The code supporting the numerical examples in this study will be made accessible upon the acceptance of the article.

\section*{Acknowledgement}
The authors were supported by the National Science Foundation of China (Nos. 12271409 and 12171367), the Natural Science Foundation of Shanghai (No. 21ZR1465800), Shanghai Municipal Science and Technology Commission (Nos. 20JC1413500, 21JC1403700 and 2021SHZDZX0100), the Interdisciplinary Project in Ocean Research of Tongji University, and the Fundamental Research Funds for the Central Universities.

\FloatBarrier

\appendix
\section{Asymptotic analysis of PI-INN}
\label{appendix-1}
Here, the asymptotic analysis of PI-INN is given in the case where all the loss terms $\mathcal{L}_{\text{equ}}$, $\mathcal{L}_{\text{bound}}$ and $\mathcal{L}_{\text{ind}}$ in Eq. (\ref{loss-function}) are zero with $N\to\infty$.

First of all, assuming $\mathcal{L}_{\mathrm{ind}} = 0$, i.e., for any $z^{'}$ and a given $[c,z]$, we have
\begin{equation}
\mathrm{log}\,q(\bm{c},z) - \mathrm{log}\,q(\bm{c},z^{'}) = \mathrm{log}\,p(z) - \mathrm{log}\,p(z^{'}),
\end{equation}
where $q(\bm c,z)$ denotes the distribution of $(\bm c,z)=g_\theta(\lambda)$ given by the INN and the prior distribution of $\lambda$.
Then,
\begin{equation}
%\begin{split}
\mathrm{log}\,\dfrac{q(z|\bm{c})q(\bm{c})}{q(z^{'}|c)q(\bm{c})} = \mathrm{log}\,\dfrac{p(z)}{p(z^{'})}\ \quad
\Rightarrow\quad\mathrm{log}\,\dfrac{q(z|\bm{c})}{p(z)} = \mathrm{log}\,\dfrac{q(z^{'}|\bm{c})}{p(z^{'})}\ .
%\end{split}
\end{equation}
Here, $\mathrm{log}\,\dfrac{q(z|\bm{c})}{p(z)}$ is a constant as $[\bm{c},z]$ is fixed. Defining $R = \mathrm{log}\,\dfrac{q(z|\bm{c})}{p(z)}$
, we have $q(z^{'}|\bm{c}) = Rp(z^{'})$, and  
\begin{equation}
\int q(z^{'}|\bm{c})\mathrm{d}z^{'} = \int p(z^{'})\mathrm{d}z^{'} = 1,\quad
\Rightarrow R = 1,\quad q(z^{'}|\bm{c}) = p(z^{'}).
\end{equation}
Therefore, $q(z|\bm{c}) \equiv p(z)$ for all $\bm{c}, z$ if $\mathcal{L}_{\mathrm{ind}} = 0$. Subsequently, we demonstrate that samples of the true posterior distribution $p(\lambda\,|\,\tilde{u})$ can be obtained using Algorithm \ref{test-algorithm}, where $\tilde{u}$ represents the observation data.

We now entertain the assumption that $\mathcal{L}_{\mathrm{ind}}$, $\mathcal{L}_{\mathrm{equ}}$, and $\mathcal{L}_{\mathrm{bound}}$ are all zero. Under this circumstance, the measurement vector $\tilde{u}$ can be exactly reconstructed using the corresponding expansion coefficient vector $\tilde{\bm{c}}$, and
\begin{eqnarray*}
p(\lambda|\tilde{u}) & = & p(\lambda|\tilde{\bm{c}})\\
 & = & \left|\mathrm{det}\left(\frac{\partial g_{\theta}(\lambda)}{\partial\lambda}\right)\right|p(\bm{c},z|\tilde{\bm{c}})\\
 & = & \left|\mathrm{det}\left(\frac{\partial g_{\theta}(\lambda)}{\partial\lambda}\right)\right|\cdot p(z)\cdot\delta(\bm{c}-\tilde{\bm{c}})\\
 & = & \left|\mathrm{det}\left(\frac{\partial g_{\theta}(\lambda)}{\partial\lambda}\right)\right|\cdot q(\bm{c},z|\tilde{\bm{c}})\\
 & = & q(\lambda|\tilde{\bm{c}})
\end{eqnarray*}
where $p(\cdot)$ denotes the distribution defined by the prior distribution
of $\lambda$ and the PDE, and $q(\cdot)$ denotes the distribution
defined by the PI-INN model. Therefore, PI-INN will
generate samples from the true posterior $p(\lambda|\tilde{u})$ by
Algorithm \ref{test-algorithm}. 

\section{Construction of the unified loss functions for data-driven and physics-informed INN}
\label{appendix-2}
In certain situations, the labeled data of $u$ corresponding to some known parameters $\lambda$ may already be available. In such cases, these pre-recorded data can be leveraged to expedite the training convergence. Here, we introduce the construction of the unified loss functions for data-driven and physics-informed INN.

Given the training sets $\{\{\lambda^{(i)},u^{(i)}\}_{i=1}^{N_{\mathrm{l}}}, 
\{\lambda^{(i)}\}_{i=1}^{N_{\mathrm{u}}}\}$, where $N_{\mathrm{l}}$ and $N_{\mathrm{u}}$ represent the number of training samples with and without labeled data, respectively. Then, the total loss function consists of two components, that is
\begin{equation}\label{general-loss-function}
\mathcal{L} = \mathcal{L}_{\mathrm{labeled}} + \mathcal{L}_{\mathrm{unlabeled}}.
\end{equation}
The first component can be written as follows
\begin{equation}\label{loss-labeled}
\mathcal{L}_{\mathrm{labeled}} = \mathcal{L}_{\mathrm{data}} + \mathcal{L}_{\mathrm{l,ind}},
\end{equation}
where
\begin{equation*}
%\begin{split}
\mathcal{L}_{\mathrm{data}} = \frac{1}{N_{\mathrm{l}}}\sum\limits_{i=1}^{N_{\mathrm{l}}}\Vert u^{(i)} - \hat{u}^{(i)}\Vert^2,\quad
\mathcal{L}_{\mathrm{l,ind}} = \dfrac{1}{N_{\mathrm{l}}}\sum\limits_{i=1}^{N_{\mathrm{l}}}\big\Vert\mathrm{log}\,q(\hat{\bm{c}}^{(i)},\hat{z}^{(i)}) - \mathrm{log}\,q(\hat{\bm{c}}^{(i)},z^{(i)}) - (\mathrm{log}\,p(\hat{z}^{(i)}) - \mathrm{log}\,p(z^{(i)}))\big\Vert^2,
%\end{split}
\end{equation*}
and the second component can be written as follows
\begin{equation}\label{loss-unlabeled}
\mathcal{L}_{\mathrm{unlabeled}} = \mathcal{L}_{\mathrm{equ}} + \mathcal{L}_{\mathrm{bound}} + \mathcal{L}_{\mathrm{u,ind}},
\end{equation}
where
\begin{equation*}
\begin{split}
\mathcal{L}_{\mathrm{equ}} &= \dfrac{1}{N_{\mathrm{u}}}\sum\limits_{i=1}^{N_{\mathrm{u}}}\big\Vert\mathcal{N}(\hat{u}^{(i)};\lambda^{(i)})\big\Vert^2,\quad
\mathcal{L}_{\mathrm{bound}} = \dfrac{1}{N_{\mathrm{u}}}\sum\limits_{i=1}^{N_{\mathrm{u}}}\big\Vert\mathcal{B}(\hat{u}^{(i)})\big\Vert^2,\\
\mathcal{L}_{\mathrm{u,ind}} &= \dfrac{1}{N_{\mathrm{u}}}\sum\limits_{i=1}^{N_{\mathrm{u}}}\big\Vert\mathrm{log}\,q(\hat{\bm{c}}^{(i)},\hat{z}^{(i)}) - \mathrm{log}\,q(\hat{\bm{c}}^{(i)},z^{(i)}) - (\mathrm{log}\,p(\hat{z}^{(i)}) - \mathrm{log}\,p(z^{(i)}))\big\Vert^2.
\end{split}
\end{equation*}
It is noted that the relative weights of loss terms are omitted for the sake of convenience. 

\section{Introduction of MMD loss term for INN}
\label{appendix-3}
The Maximum Mean Discrepancy (MMD) is a kernel-based method used to measure the distance between two probability distributions that can only be accessed through their respective samples \cite{MMD}. For instance, suppose that $X$ and $X^{'}$ are random variables following the probability distribution $p$, and that $Y$ and $Y^{'}$ are random variables following the probability distribution $q$. Then the MMD can be expressed as:
\begin{equation*}
\mathrm{MMD}[p,q] = \mathbb{E}_{X,X^{'}}[k(X,X^{'})] - 2\mathbb{E}_{X,Y}[k(X,Y)] + \mathbb{E}_{Y,Y^{'}}[k(Y,Y^{'})],
\end{equation*}
where $k$ is a kernel function. In \cite{INN}, MMD is used to measure the distribution distance between $p(\bm{c})p(z)$ and $q(\bm{c},z)$ ($\bm{c}$ refers to the measurement data in \cite{INN}. Here the notation consistent with Section \ref{networks}-\ref{training-test-approach} is used to avoid confusion). Ardizzone et al. demonstrated that if $\mathrm{MMD}[p(\bm{c})p(z), q(\bm{c},z)] = 0$, $\bm{c}$ and $z$ are independent and $z$ follows the standard Gaussian distribution. However, calculation of the MMD term requires a large amount of labeled data, which may be difficult to obtain in many practical applications\,\cite{Decreasing-Power-MMD}. 

\section{Some details implemented for numerical experiments}
\label{appendix-4}
Some settings of four numerical examples will be described in detail. For example \ref{experiment1}, ABC method and loss functions of two INN models will be introduced. For example \ref{experiment2}-\ref{experiment4}, we will focus on demonstrating the equation loss terms, and they share the same independence loss term.
\subsection{Inverse kinematics}
\label{appendix-4-1}
\subsubsection{Approximate Bayesian computation (ABC)}
\ \\
ABC is a statistical method used for approximating the likelihood function of a Bayesian model. Instead of computing the likelihood directly, ABC simulates datasets from the model and compare them to the observed data, using a distance measure to determine how well they match. Based on the principle of rejection sampling, if the distance between the simulated and observed data is smaller than a predefined threshold, the corresponding parameter values are accepted as possible parameter values for the model.

In this example, the distance threshold is $0.035$. We generate samples $\bm{x}$ from the prior distribution and select the ones where the $L_2$ distance between the corresponding simulated value $\bm{y}$ and the observation $\tilde{\bm{y}}$ is less than $0.035$ until 2000 accepted samples are obtained.

\subsubsection{Hyperparameters for INN models}
\ \\
\quad\ INN 1 consists of 8 affine coupling layers, where both $\bm{s}(\cdot)$ and $\bm{t}(\cdot)$ contain 3 fully-connected layers with 48 neurons. INN 2 has the same architecture as INN 1 for comparison (slightly different from the architecture in \cite{INN}). For two INN models, the input and output are $\bm{x} = (x_1,x_2,x_3,x_4)$ and $[y_1, y_2, z]$ with $\mathrm{ndim}(z) = 2$, respectively. The models are trained with $4000$ labeled data and a batch size of $64$. The Adam optimizer is used with an initial learning rate of $5\times 10^{-4}$ and $1200$ iterations, which decays with a rate of $0.8$ at epochs 400, 600, and 1000.

\subsection{1-d diffusion equation}
\subsubsection{Equation loss term}
Consider the following elliptic equation
\begin{equation*}
-\nabla\cdot\big(K(\bm{x})\nabla u(\bm{x})\big) = f(\bm{x}),\bm{x}\in\mathcal{D},
\end{equation*}
and the equation loss term can be constructed with its variational form:
\begin{equation}\label{galerkin-loss}
\begin{split}
\mathcal{L}_{\mathrm{equ}} &= \mathbb{E}_{K,c}\Big[\big(-\nabla\cdot(K(\bm{x})\nabla u(\bm{x})) - f(\bm{x})h(\bm{x},\bm{c},r)\big)^2\Big]\\
&= \mathbb{E}_{K,c}\Big[\big((K(\bm{x})\nabla u(\bm{x}), \nabla h(\bm{x},\bm{c},r)) - (f(\bm{x}), h(\bm{x},\bm{c},r))\big)^2\Big],
\end{split}
\end{equation}
where the inner product is defined as $(f_1, f_2) \triangleq \int_{\mathcal{D}}f_1(\bm{x})^{\top}f_2(\bm{x})\mathrm{d}\bm{x}$, and $h(\bm{x},\bm{c},r)$ is the test function:
\begin{equation*}
h(\bm{x},\bm{c},r) = r^{-D_x}1_{\Vert \bm{x}-\bm{c} \Vert_{\infty}\leq \frac{r}{2}},
\end{equation*}
where $\big\{\bm{x}|~\Vert \bm{x}-\bm{c} \Vert_{\infty}\leq \frac{r}{2}\big\}$ defines a disk in $\mathcal{D}$, $r>0$ is a constant, $c$ is randomly chosen in $\mathcal{D}$, and $D_x$ represents the spatial dimension. It can be seen that $h(\bm{x},\bm{c},r)$ defines a uniform distribution in $\big\{\bm{x}|~\Vert \bm{x}-\bm{c} \Vert_{\infty}\leq \frac{r}{2}\big\}$. How to obtain the unbiased estimate of $\mathcal{L}_{\mathrm{equ}}$\ \cite{NFFs} will be discussed as follows:

For $D_x = 1$, then
\begin{equation}
\begin{split}
\big(K(x)\nabla u(x), \nabla h(x,c,r)\big) &= r^{-1}\Big(K(c-\frac{r}{2})\nabla u(c-\frac{r}{2}) - K(c+\frac{r}{2})\nabla u(c+\frac{r}{2})\Big),\\
\big(f(x), h(x,c,r)\big) &= \mathbb{E}_{x\sim h(x,c,r)}[f(x)],
\end{split}
\end{equation}
we have

\noindent(1) For a fixed $r$, draw $c^{i}~(i=1,\cdots,n)$, thus $\{x|\Vert x-c^{i} \Vert_{\infty}\leq \frac{r}{2}\}$ is a subset of $\mathcal{D}$,

\noindent(2) Sample $x^{i}, x^{i'}\sim h(x,c^{i},r)$~$(i = 1,\cdots, n).$

\noindent(3) Given training set $\{K^{(j)}\}_{j=1}^{N}$, calculate 
\begin{equation}
\big(K(x)\nabla u(x), \nabla h(x,c^{i},r)\big) = r^{-1}\Big(K(c^{i}-\frac{r}{2})\nabla u(c^{i}-\frac{r}{2}) - K(c^{i}+\frac{r}{2})\nabla u(c^{i}+\frac{r}{2})\Big),
\end{equation}
for $j = 1,\cdots,N, i = 1, \cdots n.$

\noindent(4) Calculate the unbiased estimate of $\mathcal{L}_{\mathrm{equ}}$:
\begin{equation}
\hat{\mathcal{L}}_{\mathrm{equ}} = \frac{1}{nN}\sum\limits_{j=1}^{N}\sum\limits_{i=1}^{n}\mathcal{L}_{\mathrm{equ}}^{(j),i},
\end{equation}
where
\begin{equation}
\mathcal{L}_{\mathrm{equ}}^{(j),i} = \big(\big(K^{(j)}(x)\nabla u(x),\nabla h(x,c^{i},r)\big) - f(x^{i})\big)\cdot\big(\big(K^{(j)}(x)\nabla u(x),\nabla h(x,c^{i},r)\big) - f(x^{i'})\big).
\end{equation}

We set $u(x;\lambda)=\sum_{i=1}^{5}c_i(\lambda)\phi_{i}(x)$
Assumed that $z$ follows a 5-dimensional standard Gaussian distribution, and the spatial coordinates $\{\bm{x}^{(j)}\}_{j=1}^{M}$ for training are the same as the grid points for the FEM.

\subsubsection{Importance sampling technique}
\label{appendix-IP}
\ \\
In our study, we enhance sampling accuracy by employing the importance sampling (IP) method, leveraging the informative proposal distribution generated by PI-INN. The implementation steps are as follows:
\begin{itemize}
    \item 1. Sampling according to Algorithm \ref{test-algorithm}: $\{\tilde{\lambda}^{i}\}_{i=1}^{N}$;
    \item 2. Calculate the weights: $w_i = p(\tilde{\lambda}^{(i)}\,|\,\tilde{u}) / q(\tilde{\lambda}^{(i)}\,|\,\tilde{u}), i=1,\cdots,N$, where $p(\tilde{\lambda}^{(i)}\,|\,\tilde{u})$ and $q(\tilde{\lambda}^{(i)}\,|\,\tilde{u})$ are densities of the exact Bayesian posterior distribution and the proposal distribution by PI-INN, respectively.
    \item 3. Normalize the weights to ensure they sum to 1: $w_i\leftarrow w_i / \sum\limits_{i=1}^{N}w_i$;
    \item 4. Calculate the statistics: $\mathbb{E}_{p(\lambda|\tilde{u})}(f(\tilde{\lambda})) \approx \sum\limits_{i=1}^{N}w_if(\tilde{\lambda}^{i})\hat{=}\tilde{\lambda}_{IP}$.
\end{itemize}
% We would like to emphasize that the low sampling cost associated with drawing samples from $q(\tilde{\lambda},|,\tilde{u})$ permits us to iterate the aforementioned steps multiple times, yielding a sample set for $p(\tilde{\lambda}^{(i)}\,|\,\tilde{u})$: $\{\tilde{\lambda}_{IP}^{(k)}\}_{k=1}^{N_{\mathrm{sample}}}$.

\subsubsection{Hyperparameters for numerical examples}
\begin{itemize}
    \item GRF prior
\end{itemize}
The INN consists of 8 coupling layers, where both $\bm{s}(\cdot)$ and $\bm{t}(\cdot)$ have 1 hidden layer with 100 neurons and ReLU activation function. The NB-Net has 6 hidden layers and 64 neurons.
We use $4000$ unlabeled data with a batch size of $64$ for training and train our model for 1000 epochs using the Adam optimizer in Pytorch. The learning rate is initially set to $1\times 10^{-3}$ and decays with a rate of 0.8 at epochs 600 and 900.

\begin{itemize}
    \item Non-Gaussian mixed prior
\end{itemize}
We use the same hyperparameters for the architecture as the GRF prior case. Utilizing the Adam optimizer \cite{Adam} with an initial learning rate of 0.001, the PI-INN model underwent training for $2400$ iterations, with a learning rate decay rate of $0.8$ applied every $400$ iterations.

\subsubsection{Proposal distributions for MCMC}
\ \\
GRF prior \& Mixed Non-Gaussian prior: Gaussian distribution with zero mean and standard deviation of 0.01;

\subsubsection{Definition of the absolute mean point-wise error}
\label{appendix-4-pointwise-error}

\ \\
\begin{equation}
\mathrm{Mean\ pointwise\ error} = \sum\limits_{k=1}^{K}\vert f(\{\lambda^{(i)}\}_{i=1}^{N})(x_k) - f^{*}(x_k)\vert / K,
\end{equation}
where $\{\lambda^{(i)}\}_{i=1}^{N}$ is the sample set for posterior distribution, $f(\{\lambda^{(i)}\}_{i=1}^{N})$ represents a certain statistic for posterior samples, $f^{*}$ is the exact value for $f$ with exact posterior samples.

\subsubsection{Effect of noise intensity for Bayesian inference}
\label{appendix-4-2-2}
\begin{figure}[htp]
	\centering
	\includegraphics[scale=0.32]{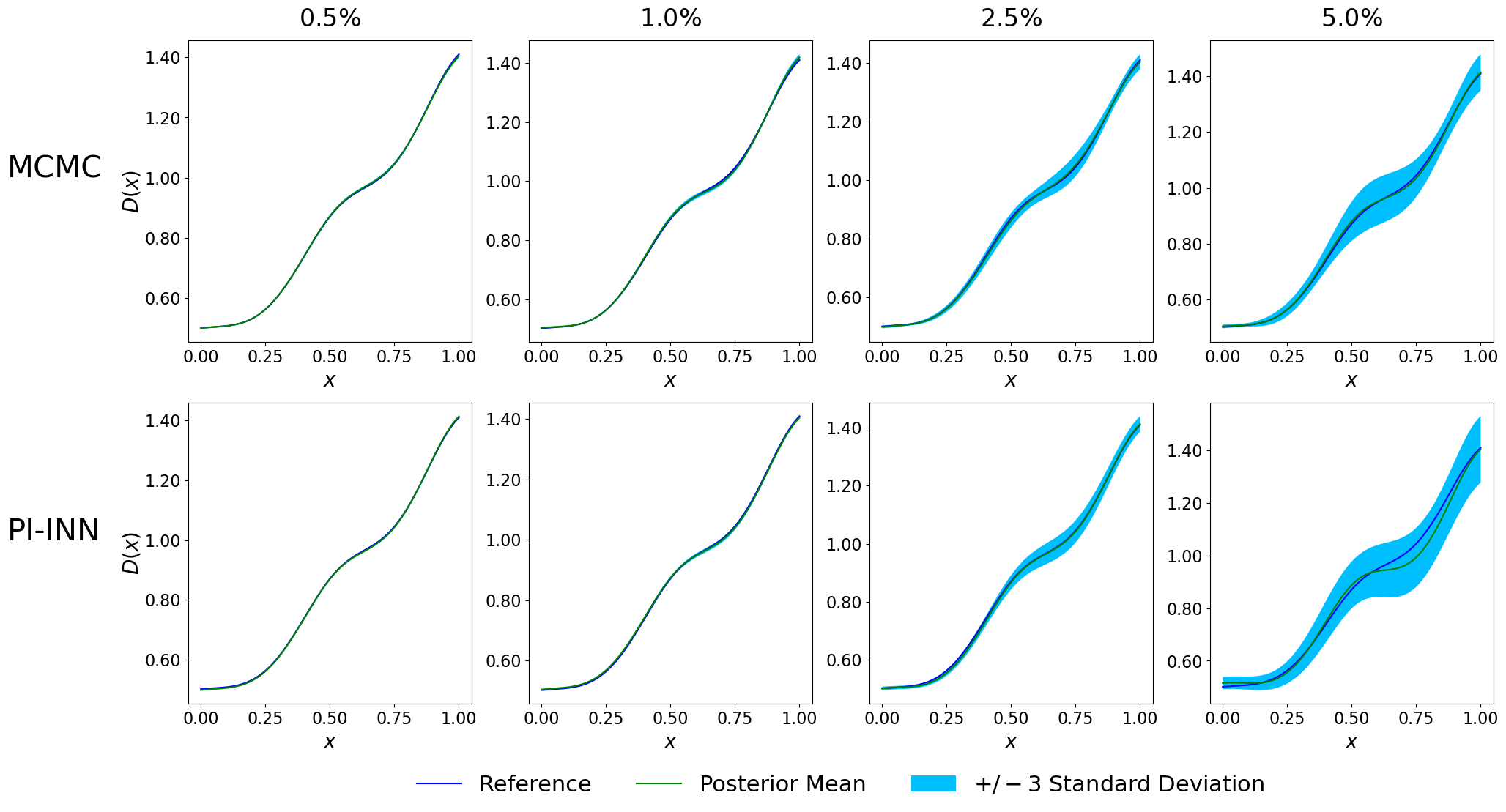}
	\caption{Comparison of posterior samples by PI-INN and reference MCMC method with increasing observation noise. The first and second rows correspond to the MCMC and PI-INN results, respectively. The columns represent the posterior samples under observation noise with different standard deviations, specifically, $0.005$, $0.01$, $0.025$ and $0.05$, respectively.}
	\label{noise-effect}
\end{figure}

\ \\
As illustrated in Fig.\,\ref{ex2-diffusion-guassian-slices} and Fig.\,\ref{ex2-diffusion-mixed-guassian-slices}, our model consistently yields posterior samples that align closely with those obtained through MCMC when zero-mean Gaussian noise with $0.01$ standard deviation is added to the observation data. However, it is imperative to acknowledge that our current method is suitable for the low-noise scenarios characterized by low noise levels. To substantiate this assertion, we examine the impact of increasing noise intensity on the posterior samples, as displayed in Fig.\,\ref{noise-effect}. In this analysis, we adopt the exponential form of GRF as the prior for $D(x)$ in Eq.\,(\ref{ex2-equ}). We observe that PI-INN successfully captures the actual posterior distribution when the standard deviation of Gaussian noise remains below $0.025$, which is consistent with the outcomes depicted in Fig. \ref{ex2-noise-sensor-analysis}. Nevertheless, as the standard deviation increases to $0.05$, the accuracy of UQ for $D(x)$ diminishes.

\subsection{2-d Darcy equation}
\subsubsection{Equation loss term}
\ \\
The equation loss term similar to Eq. (\ref{galerkin-loss}) is also used for the Darcy equation. For $D_x = 2$, let $\bm{x} = (x_1,x_2)^{\top}$ and $\bm{c} = (c_1,c_2)^{\top}$, we have
\begin{equation}
\nabla h(x_1,x_2,\bm{c},r) = \begin{pmatrix}
r^{-2}\big(\delta(x_1 - c_1 + \frac{r}{2}) - \delta(x_1 - c_1 - \frac{r}{2})\big)1_{\vert x_2 - c_2\vert\leq\frac{r}{2}}\\
r^{-2}\big(\delta(x_2 - c_2 + \frac{r}{2}) - \delta(x_2 - c_2 - \frac{r}{2})\big)1_{\vert x_1 - c_1\vert\leq\frac{r}{2}},
\end{pmatrix}
\end{equation}
Then,
\begin{equation}
\begin{split}
&\big(K(x_1,x_2)\nabla u(x_1,x_2), \nabla h(x_1,x_2,\bm{c},r)\big)\\
&= \int\int K(x_1,x_2)\nabla_{x_1}u(x_1,x_2)\dfrac{\delta(x_1-c_1+\frac{r}{2}) - \delta(x_1-c_1-\frac{r}{2})}{r^2}1_{\vert x_2-c_2\vert\leq\frac{r}{2}}\mathrm{d}x_1\mathrm{d}x_2\\
&+ \int\int K(x_1,x_2)\nabla_{x_2}u(x_1,x_2)\dfrac{\delta(x_2-c_2+\frac{r}{2}) - \delta(x_2-c_2-\frac{r}{2})}{r^2}1_{\vert x_1-c_1\vert\leq\frac{r}{2}}\mathrm{d}x_1\mathrm{d}x_2\\
&= \mathbb{E}_{x_2\sim\mathcal{U}(c_2-\frac{r}{2},c_2+\frac{r}{2})}\bigg[\dfrac{K(c_1-\frac{r}{2},x_2)\nabla_{x_1}u(c_1-\frac{r}{2},x_2) - K(c_1+\frac{r}{2},x_2)\nabla_{x_1}u(c_1+\frac{r}{2},x_2)}{r}\bigg]\\
&+\mathbb{E}_{x_1\sim\mathcal{U}(c_1-\frac{r}{2},c_1+\frac{r}{2})}\bigg[\dfrac{K(x_1,c_2-\frac{r}{2})\nabla_{x_2}u(x_1,c_2-\frac{r}{2}) - K(x_1,c_2+\frac{r}{2})\nabla_{x_2}u(x_1,c_2+\frac{r}{2})}{r}\bigg].
\end{split}
\end{equation}
Then, the analysis of the unbiased estimate of $\mathcal{L}_{\mathrm{equ}}$ can be given as follows:

\noindent(1) For fixed $r$, sample $\bm{c}^{i} = (c_1^{i},c_2^{i})$ so that $\{\bm{x}|\Vert \bm{x}-\bm{c}^{i} \Vert_{\infty}\leq \frac{r}{2}\}\in\mathcal{D}$

\noindent(2) Draw $x_1^{i}, x_1^{i'}\sim\mathcal{U}(c_1^{i}-\dfrac{r}{2},c_1^{i}+\dfrac{r}{2})$ and $x_2^{i}, x_2^{i'}\sim\mathcal{U}(c_2^{i}-\dfrac{r}{2},c_2^{i}+\dfrac{r}{2})$ for $i = 1,\cdots,n$

\noindent(3) Given training set $\{K^{(j)}\}_{j=1}^{N}$, calculate
\begin{equation}
\begin{split}
e^{(j),i} &= r^{-1}K^{(j)}(c_1^{i}-\frac{r}{2},x_2^{i})\nabla_{x_1}u(c_1^{i}-\frac{r}{2},x_2^{i}) - r^{-1}K^{(j)}(c_1^{i}+\frac{r}{2},x_2^{i})\nabla_{x_1}u(c_1^{i}+\frac{r}{2},x_2^{i})\\
&+ r^{-1}K^{(j)}(x_1^{i},c_2^{i}-\frac{r}{2})\nabla_{x_2}u(x_1^{i},c_2^{i}-\frac{r}{2}) - r^{-1}K^{(j)}(x_1^{i},c_2^{i}+\frac{r}{2})\nabla_{x_2}u(x_1^{i},c_2^{i}+\frac{r}{2})\\
&- f(x_1^{i},x_2^{i})
\end{split}
\end{equation}
\begin{equation}
\begin{split}
e^{(j),i'} &= r^{-1}K^{(j)}(c_1^{i'}-\frac{r}{2},x_2^{i'})\nabla_{x_1}u(c_1^{i'}-\frac{r}{2},x_2^{i'}) - r^{-1}K^{(j)}(c_1^{i'}+\frac{r}{2},x_2^{i'})\nabla_{x_1}u(c_1^{i'}+\frac{r}{2},x_2^{i'})\\
&+ r^{-1}K^{(j)}(x_1^{i'},c_2^{i'}-\frac{r}{2})\nabla_{x_2}u(x_1^{i'},c_2^{i'}-\frac{r}{2}) - r^{-1}K^{(j)}(x_1^{i'},c_2^{i'}+\frac{r}{2})\nabla_{x_2}u(x_1^{i'},c_2^{i'}+\frac{r}{2})\\
&- f(x_1^{i'},x_2^{i'})
\end{split}
\end{equation}
for $j = 1,\cdots,N, i = 1,\cdots,n$.

\noindent(4) Calculate the unbiased estimate of $\mathcal{L}_{\mathrm{equ}}$:
\begin{equation}
\hat{\mathcal{L}}_{\mathrm{equ}} = \frac{1}{nN}\sum\limits_{j=1}^{N}\sum\limits_{i=1}^{n}e^{(j),i}e^{(j),i'}.
\end{equation}

\subsubsection{Hyperparameters for numerical example}
\ \\
The INN consists of 8 coupling layers, where both $\bm{s}(\cdot)$ and $\bm{t}(\cdot)$ contain 1 hidden layer with 100 neurons and ReLU activation functions. The NB-Net is constructed by a fully-connected neural network with $5$ hidden layers containing $128$ neurons and ReLU activation function. The spatial coordinates $\{\bm{x}^{(j)}\}_{j=1}^{M}$ for training are consistent with the grid points used for FEM. The variational formulation of $\mathcal{L}_{\mathrm{equ}}$ is still adopted. $8000$ unlabeled data with a batch size of $64$ is utilized. The training of the model was conducted using the Adam optimizer, initially set with a learning rate of $5\times 10^{-4}$, which decayed by a factor of $0.8$ every $800$ iterations. The training process spanned $5000$ iterations in total.

\subsection{2-d seismic tomography}
\label{appendix-4-4}
\subsubsection{Equation loss term}
\ \\
Considering the singularity of the single point source and substituting $T(\bm{x}) = T_0(\bm{x})\tau(\bm{x})$ into Eq. (\ref{eikonal-equ}), the following factored eikonal equation will be used \cite{Factored-eikonal}:
\begin{equation}\label{factor-eikonal}
\begin{array}{cc}
T_0^2\vert\nabla\tau\vert^2 + \tau^2\vert\nabla T_0\vert^2 + 2T_0\tau(\nabla T_0\nabla\tau) = \dfrac{1}{v^2(\bm{x})},\quad 
\tau(\bm{x}_s) = 1,
\end{array}
\end{equation}
where $T_0(\bm{x}) = \dfrac{\vert\bm{x} - \bm{x}_s\vert}{v(\bm{x}_s)}$. Then, the equation loss can be written as the following two terms:

\begin{equation}
\mathcal{L}_{\mathrm{equ}} = \mathcal{L}_{\mathrm{factored}} + \mathcal{L}_{\mathrm{source}},
\end{equation}
where
\begin{equation*}
\begin{split}
\mathcal{L}_{\mathrm{factored}} &= \dfrac{1}{N}\sum\limits_{i=1}^{N}\big\Vert T_0^2\vert\nabla\hat{\tau}^{(i)}\vert^2 +(\hat{\tau}^{(i)})^2\vert\nabla T_0\vert^2 + 2T_0\hat{\tau}^{(i)}(\nabla T_0\nabla\hat{\tau}^{(i)}) - 1/(v^{(i)})^2\big\Vert^2,\\
\mathcal{L}_{\mathrm{source}} &= \dfrac{1}{N}\sum\limits_{i=1}^{N}\Vert\hat{\tau}^{(i)} - 1\Vert^2,
\end{split}
\end{equation*}
where $\{v^{(i)}\}_{i=1}^{N}$ is the training set, $\hat{\tau}^{(i)}$ is the prediction solutions of $\tau^{(i)}$ by PI-INN. It is noted that $\mathcal{L}_{\mathrm{bound}}$ is not used for this example. 

\subsubsection{Hyperparameters for numerical example}
For the network architecture, the INN comprises 8 coupling layers, where both $\bm{s}(\cdot)$ and $\bm{t}(\cdot)$ have 1 hidden layer with 100 neurons. The NB-Net contains $6$ hidden layers with $128$ neurons, and the ReLU activation function is utilized for both INN and the NB-Net. Using the Adam optimizer with a learning rate of $5\times 10^{-4}$, we train the model over $N = 4000$ unlabeled data with a batch size of $64$.

% For the dimensionality of random variables, we set $\mathrm{ndim}(\bm{c}) = 4$ and $\mathrm{ndim}(z) = 6$. 
% Since seven dimensional random parameters are considered for velocity models, padding variables $\zeta$ with $\mathrm{ndim}(\zeta) = 3$ have been added.

\newpage
\section*{References}
\bibliography{refs}

\end{document}